\documentclass{article}
\usepackage[T1]{fontenc}
\usepackage{a4}
\usepackage[dvips]{graphicx}
\usepackage{amsmath}
\usepackage{amsthm}
\usepackage{amssymb}
\usepackage{amsfonts}
\usepackage{mathrsfs}
\usepackage[all]{xy}
\newtheorem{etheo}[subsubsection]{Theorem}
\newtheorem{ecor}[subsubsection]{Corollary}
\newtheorem{econj}[subsubsection]{Conjecture}
\newtheorem{elem}[subsubsection]{Lemma}
\newtheorem{eprop}[subsubsection]{Proposition}
\newtheorem{etheosec}[subsection]{Theorem}

\newtheorem{econjsec}[subsection]{Conjecture}

\newtheorem{edefin}[subsubsection]{Definition}

\newtheorem{eexm}[subsubsection]{Example}
\newtheorem{erem}[subsubsection]{Remark}

\newcommand{\A}{\mathbb{A}}
\newcommand{\C}{\mathbb{C}}
\newcommand{\R}{\mathbb{R}}
\newcommand{\Q}{{\mathbb{Q}}}

\newcommand{\qb}{{\overline{\Q}}}
\newcommand{\Z}{\mathbb{Z}}
\newcommand{\N}{\mathbb{N}}
\newcommand{\p}{\mathfrak{p}}
\newcommand{\q}{\mathfrak{q}}

\newcommand{\U}{\mathbb{U}}

\newcommand{\G}{\mathbb{G}}
\newcommand{\qlb}{{\overline{\Q_\ell}}}
\newcommand{\ql}{{\Q_\ell}}
\newcommand{\qp}{{\Q_p}}
\newcommand{\zp}{{\Z_p}}
\newcommand{\zl}{{\Z_\ell}}
\newcommand{\II}{\mathscr{I}}
\newcommand{\DD}{\mathscr{D}}

\newcommand{\GKK}{\mathscr{G}_K}
\newcommand{\GLL}{\mathscr{G}_L}

\newcommand{\Spec}{\mathrm{Spec}}

\newcommand{\id}{\mathrm{id}}
\newcommand{\Aut}{\mathrm{Aut}}
\newcommand{\Hom}{\mathrm{Hom}}

\newcommand{\End}{\mathrm{End}}
\newcommand{\norm}{\mathrm{norm}}

\newcommand{\Gal}{\mathrm{Gal}}
\newcommand{\Card}{\mathrm{Card}}
\newcommand{\ad}{\mathrm{ad}}
\newcommand{\GL}{\mathrm{GL}}

\newcommand{\Lie}{\mathrm{Lie}}
\newcommand{\Gr}{\mathrm{Gr}}

\newcommand{\GSp}{\mathrm{GSp}}
\newcommand{\MT}{\mathrm{MT}}
\newcommand{\Sh}{\mathrm{Sh}}
\newcommand{\limproj}{\underset{\gets}{\lim}\,}
\newcommand{\Hun}{\mathrm{H_{1}}}
\newcommand{\hunah}{\mathrm{h_{1,AH}}}
\newcommand{\Rep}{\underline{\mathrm{Rep}}}
\newcommand{\Res}{\mathrm{Res}}
\newcommand{\MF}{\underline{\mathrm{MF}}}
\newcommand{\Vect}{\underline{\mathrm{Vect}}}

\newcommand{\D}{\underline{\mathrm{D}}}
\newcommand{\V}{\underline{\mathrm{V}}}
\newcommand{\PEL}{\textrm{PEL}}

\newcommand{\Fil}{\mathrm{Fil}}

\title{Galois representations, Mumford-Tate groups and good reduction of abelian varieties}

\author{Fr\'ed\'eric Paugam}

\begin{document}
\maketitle
\abstract{Let $K$ be a number field and $A$ an abelian variety over $K$. We are interested in the following conjecture of Morita: if the Mumford-Tate group of $A$ does not contain unipotent $\Q$-rational points then $A$ has potentially good reduction at any discrete place of $K$. The Mumford-Tate group is an object of analytical nature whereas having good reduction is an arithmetical notion, linked to the ramification of Galois representations.
This conjecture has been proved by Morita for particular abelian varieties with
many endomorphisms (called of {\PEL} type). Noot obtained results for abelian
varieties without nontrivial endomorphisms (Mumford's example, not of {\PEL} type). We give new results for abelian varieties not of {\PEL} type.}



\section*{Introduction}
Let $E$ be a number field, $\p$ a finite prime of $E$, $\mathcal{O}_{E,\p}$ the localization of the ring of integers at $\p$. Fix an embedding $\iota$ of $E$ in $\C$. Let $A/E$ be an abelian variety. We will denote $\Hun(A,\Q)$ the rational singular homology group of the complex abelian variety $A_{\C,\iota}$ obtained by scalar extension through $\iota$.

We will say that \emph{$A$ has good reduction at $\p$} if there exists an abelian scheme $\mathcal{A}$ over $\mathcal{O}_{E,\p}$ such that $\mathcal{A}\times_{\Spec(\mathcal{O}_{E,\p})} \Spec(E)\cong A$. We will say that $A$ has potentially good reduction at $\p$ if there exists a finite extension $E'$ of $E$ and a prime $\p'$ of $E'$ over $\p$ such that $A_{E'}$ has good reduction at $\p'$.

Let $\mathbb{S}=\Res_{\C/\R}\G_{m,\C}$ be the restriction of scalars from $\C$ to $\R$ of the multiplicative group. Define the \emph{weight cocharacter} $w\colon \G_{m,\R}\to \mathbb{S}$ to be the morphism given on points by the naturel inclusion $\R^*=\G_{m,\R}(\R)\subset\mathbb{S}(\R)=\C^*$.
A \emph{$\Q$-Hodge structure} is a $\Q$-vector space $V$ endowed with a morphism $h\colon \mathbb{S}\to \GL(V_\R)$ such that $h\circ w\colon \G_{m,\R}\to \GL(V_\R)$ is defined over $\Q$.
The \emph{Mumford-Tate group of a Hodge structure $(V,h)$} is the smallest $\Q$-algebraic group $\MT(V,h)\subset~\GL(V)$ that contains the image of $h$ after scalar extension to $\R$.
Hodge theory tells us that the singular homology $\mathrm{H}_i(X,\Q)$ of every projective smooth algebraic variety $X$ over $\C$ carries a $\Q$-Hodge structure. In the case of the abelian variety $A$, the Hodge structure on $\Hun(A,\Q)$ is given by the complex structure on the real torus
$$
A_\C\cong \Hun(A,\R)/\Hun(A,\Z).
$$
The Mumford-Tate group of $A$ is defined as the Mumford-Tate group of the Hodge structure $\Hun(A,\Q)$.

Our work deals with the following conjecture, first stated by Morita in \cite{Morita}, p437.
\begin{econjsec}[Morita]
\label{conjalgmorita}
If the Mumford-Tate group of $A$ contains no nontrivial unipotent $\Q$-rational point then $A$ has potentially good reduction at every finite prime of $E$.
\end{econjsec}

The original conjecture was formulated in more geometrical terms that we will explain without going into details of the theory of Shimura varieties, referring to section \ref{constab} for more precisions on this.

The couple $(\MT(A),h\colon \mathbb{S}\to \MT(A)_\R)$ is a Shimura datum, that will permit to construct (the canonical model of) a Shimura variety $\Sh$ over some number field. We can assume that this variety is defined over a finite extension $E'$ of $E$. To the natural representation $\MT(A)\to \GL(\mathrm{H}_1(A_\C^{an},\Q))$, we can associate an abelian scheme $\mathcal{A}$ over $\Sh$ giving a family of abelian varieties which ``contains $A_{E'}$'' in the following sense: there exists $x\in\Sh(E')$ and a cartesian diagram
$$
\xymatrix{
A_{E'}\ar[r]\ar[d]     	& \mathcal{A}\ar[d]\\
\Spec(E')\ar[r]^{x}	& \Sh
}.
$$

The variety $\Sh$ is a moduli space for polarized abelian varieties endowed with a level structure and an (interesting) \emph{additional structure}. Roughly speaking, it is a moduli space for abelian varieties with Mumford-Tate group contained in $\MT(A)$.

It will be shown in \ref{lem-compacite} that the hypothesis on unipotents in Morita's conjecture is equivalent to the compactness of $\Sh_\C$.

In the {\PEL} case, treated partially by Morita, the additional structure that defines the moduli problem is given by endomorphisms of the abelian variety. So the Shimura variety is a moduli space for abelian varieties endowed with a \textbf{P}olarization, an action of some \textbf{E}ndomorphism ring, and a \textbf{L}evel structure. This kind of moduli problem can be defined in finite characteristic, which makes the conjecture easier to handle.

Noot \cite{Noot2} proved the conjecture in the case of Mumford's Shimura curves, that are non {\PEL} Shimura varieties. In this case, the additional structure used to define the moduli problem for $\Sh$ is not given by endomorphisms of the abelian variety but by (co)homology cycles (more precisely absolute Hodge cycles). A difficulty with this kind of moduli problem is that it is not easy to define in finite characteristic because absolute Hodge cycles are not easy to reduce modulo some prime (unless one admits the Hodge conjecture).

We now sketch the results we obtain with respect to Morita's conjecture. Take an abelian variety $A$ defined over a number field $E\subset \C$. Let $V=\Hun(A_\C,\Q)$ be its first singular homology group and $\rho\colon \MT(A)\to \GL(V)$ be the natural representation of its Mumford-Tate group.

We show in a first section \ref{theo-criterion} that if the image of the representation $\rho_\ql$ obtained by scalar extension from $\Q$ to $\ql$, contains no nontrivial unipotent element over $\ql$ of index $2$, then the abelian variety has potentially good reduction at every discrete place of its field of definition. This result is based on a theorem of Deligne about absolute Hodge cycles and on our study in \ref{expo-monod-padic} of an etale $p$-adic analog of the exponential of the monodromy in $\ell$-adic cohomology.

In a second part, we classify the representations satisfying the hypothesis of our theorical criterion in terms of the combinatorics of the Galois action on their tensorial and irreducible components. These combinatorics are encoded in a combinatorial datum called a polymer. We thus obtain a new criterion of good reduction based on combinatorial properties of the representation of the Mumford-Tate group.

In a third part, we use the machinery of Shimura varieties to construct families of abelian varieties satisfying our criterion.
We also show that our results are different from Morita's results because the Shimura varieties we consider are not of {\PEL} type.

In a fourth part, we give a transposition method that permits to apply our results to abelian varieties with same adjoint Mumford-Tate groups as abelian varieties satisfying our criterion. We also give examples of abelian varieties for which this transposition method works.

\emph{Acknowledgements.} This article is a part of my doctoral thesis.
I would like to thank my Phd professor Rutger Noot for his help and
for proposing me to work on this exciting subject. I also thanks
J.-P. Wintenberger, Ben Moonen and Yves Andr\'e for useful discussions and
comments. I would also like to thank the referee for useful simplification
suggestions. To finish, I am grateful to Rennes' s laboratory for excellent
working conditions.

\section{Theoretical criterion for potentially good reduction}
\label{theo-criterion}
Let $K$ be a field of characteristic $0$, complete for a discrete valuation $v$,  of perfect residue field $k$ of characteristic $p$. Let $K_0$ be the fraction field of the ring of Witt vectors with coefficients in $k$, $\sigma$ be the absolute Frobenius of $K_0$ (i.e. the unique continuous automorphism inducing $x\mapsto x^p$ on $k$) and $\GKK=\Gal(\overline{K}/K)$. Let $\bar{v}$ be a place of $\overline{K}$ over $v$ and $\II(\bar{v})$ be the inertia subgroup at $\bar{v}$ of $\GKK$.

Let $\ell$ be a prime, $A$ be an abelian variety over $K$ and let $A[\ell^n]$ denote the $\ell^n$-torsion points in $A(\overline{K})$. The $\ell$-adic rational Tate module of $A$ is the $\ql$-vector-space $V_\ell(A)=(\limproj A[\ell^n])\otimes_{\zl}\ql$.

For $\ell\neq p$, the Neron-Ogg-Shafarevich criterion (see \cite{SeTa}, Th 1) permits to translate the question of good reduction of $A$ at $v$ in terms of the Galois representation $V_\ell(A)$.
\begin{etheosec}
The abelian variety $A$ has good reduction at $v$ if and only if the inertia subgroup $\II(\bar{v})$ of $\GKK$ acts trivially on $V_\ell(A)$.
\end{etheosec}

The aim of this paragraph is to construct a unipotent subgroup of $\GL(V_p(A))$
that contains only unipotents of index $2$ (i.e. a subgroup whose elements verify the equality $(g-\id)^2=0$), that is contained in the Zariski closure of the image of the $p$-adic Galois representation $V_p(A)$ and that will play the role of the inertia subgroup in the last theorem in the case $\ell=p$. It will be called the \emph{$p$-adic exponential of the monodromy}.
To construct it, we will work with semi-stable representations of $\GKK$ and with their (log-)cristalline counterparts:  Fontaine's filtered $(\Phi,N)$-modules.
 
Roughtly speaking, our construction gives a $p$-adic counterpart to the construction of the $\ell$-adic exponential of the monodromy \cite{Ill2}, 1.5.2 that we will briefly explain geometrically now.

For $\ell\neq p$, if $A$ is semi-stable over $K$, the action of the inertia subgroup $\II(\bar{v})$ on $V_\ell(A)$ is unipotent and is given by the exponential of a morphism
$$N\colon V_\ell(A)(1)\to V_\ell(A)$$
of $\ell$-adic representations of $\GKK$ (see \cite{Ill2}, 1.5.2 for details).

There is an analog of this morphism on the log-cristalline cohomology of the special fibre of a proper model of $A$ over $\mathcal{O}_v$ with semi-stable reduction. This gives the $N$ in the filtered $(\Phi,N)$-module of $A$.

We will prove in this section that for abelian varieties, this morphism in the log-cristalline setting gives a morphism
$$N\colon V_p(A)(1)\to V_p(A)$$
of $p$-adic representations of the local Galois group $\GKK$ that is called the \emph{$p$-adic monodromy morphism}. We must add that there is no hope to have such a $p$-adic monodromy morphism for every semi-stable $p$-adic Galois representation. The transversality of monodromy is studied in \cite{Breuil}
and there is a counter example with the galois representation of a modular
form given in the review of this article on mathscinet.

Our proof is based on a simple property of the $(\Phi,N)$-module of $A$ analogous to Griffiths transversality and this allows us to make a construction that is functorial and compatible with tensor-product and usual tannakian operations.

The first example of a semi-stable representation is $\qp(1)$, the $p$-adic galois representation given by the Galois action on $p^n$ roots of unity in $\overline{K}$ for all $n$: 
$$\qp(1)=\qp\otimes_{\Z_p}(\limproj \mathbb{\mu}_{p^n}(\overline{K})).$$
This is a one-dimensional vector-space over $\qp$.

Let $\qp(1)_0$ be the vector space $\qp$ endowed with the action of $\GKK$ given by the cyclotomic character $\chi\colon \GKK\to\zp^*$.
We will choose once for all an isomorphism $\psi\colon \qp(1)\to\qp(1)_0$. This isomorphism is determined by the element $t_{\psi}=\psi^{-1}(1)$ of $\qp(1)$ and is called a \emph{trivialisation of the Tate twist}.

\subsection{Fontaine's modules}
We recall here the definitions (see Fontaine \cite{Fontaine}, 4.3.2). 
Notations are as above.

\begin{edefin}
We call \emph{filtered $(\Phi,N)$-module} the data of a finite dimensional $K_0$-vector space $D$ endowed with
\begin{enumerate}
\item an injective $\sigma$-semi-linear map (Frobenius)
$$\Phi\colon D\to D,$$
\item a $K_0$-linear endomorphism (monodromy)
$$N\colon D\to D,$$
\item a filtration $(\Fil^iD_K)_{i\in \Z}$ of $D_K=D\otimes_{K_0}K$ by sub-$K$-vector spaces, that is required to be decreasing ($\Fil^{i+1}D_K\subset \Fil^iD_K$), exhaustive ($\cup \Fil^iD_K=D_K$) and separated ($\cap \Fil^iD_K=0$).
\end{enumerate}
We impose on these data the following compatibility condition: 
\begin{itemize}
\item $N\Phi=p\Phi N$.
\end{itemize}
We denote by $\MF_K(\Phi,N)$ the category of such modules.
\end{edefin}

The power of Fontaine's theory is that it permits to reduce problems about $p$-adic Galois representations coming from algebraic geometry, i.e. from the $p$-adic etale cohomology of projective smooth algebraic varieties, to problems of semi-linear algebra. We now define a category of $p$-adic Galois representations that includes the representations coming from algebraic geometry.

Let $V$ be a continuous finite dimensional $p$-adic representation of $\GKK=\Gal(\overline{K}/K)$. By definition, this is a continuous representation of $\GKK$ in a finite dimensional $\qp$-vector-space. We denote
$$\D_{st,K}(V)=(B_{st,K}\otimes_{\qp}V)^{\GKK},$$
where $B_{st,K}$ denotes Fontaine's period ring, defined in \cite{Fontaine2}.

We have $\dim_{\qp}(V)\geq \dim_{K_0}(\D_{st,K}(V))$ and
$V$ is called semi-stable if
$$\dim_{\qp}(V)=\dim_{K_0}(\D_{st,K}(V))$$
and potentially semi-stable if there is a finite extension $L/K$ such that
$$\dim_{\qp}(V)=\dim_{L_0}(\D_{st,L}(V)).$$
In the last equality, $V$ denotes the $p$-adic representation obtained by restriction to $\Gal(\overline{K}/L)$ and $L_0$ is the fraction field of the ring of Witt vectors of the residue field of $L$.
We denote by $\Rep_{st}(\GKK)$ (resp. $\Rep_{pst}(\GKK)$) the category of semi-stable (resp. potentially semi-stable) representations of $\GKK$. Those are tannakian categories over $\qp$ endowed with a natural fibre functor $\omega_p\colon \Rep_{pst}(\GKK)\to \Vect_\qp$ which associates to any $p$-adic representation the underlying $\qp$-vector space. The $C_{pst}$ conjecture, which is now a theorem (see \cite{Berthelot2}, 6.3.3 and references therein) claims that every $p$-adic Galois representation coming from the geometry of projective smooth varieties is potentially semi-stable.

The essential image of $\D_{st,K}$ is denoted by $\MF_K^a(\Phi,N)$. Its objects are called admissible filtered $(\Phi,N)$-modules. This is also a tannakian category over $\qp$ and the functor
$$\D_{st}\colon \Rep_{st}(\GKK) \to \MF_K^a(\Phi,N)$$
is a tensor equivalence of tannakian categories (See \cite{Fontaine}, 5.6.7, 5.6.8). A quasi-inverse $\V_{st}$ is given by $\V_{st}(D):=\Hom_{\MF_K^a(\Phi,N)}(K_0,B_{st,K}\otimes D)$. There is a natural fibre functor $\omega_{st}\colon \MF_K^a(\Phi,N)\to \Vect_{K_0}$ which associates to any filtered $(\Phi,N)$-module the underlying $K_0$-vector-space.

\begin{eexm}
\label{exmtate}
We denote by $K(1)$ the filtered $(\Phi,N)$-module $\D_{st,K}(\qp(1))$. The trivialisation $\psi\colon \qp(1)\to \qp(1)_0$ of the Tate twist corresponds to the choice of a basis $t_\psi=\psi^{-1}(1)$ for $\qp(1)$ and gives a basis $t_\psi^{-1}\otimes t_\psi$ of $\D_{st,K}(\qp(1))=(B_{st,K}\otimes_{\qp}\qp(1))^{\GKK}$ and an isomorphism $i_\psi\colon K(1)\to K(1)_0$ with the module $K(1)_0$ whose underlying vector space is $K_0$, endowed with the Frobenius $\Phi=(1/p)\sigma$, the trivial monodromy $N=0$ and the filtration $\Fil^iK(1)_0=K(1)_0$ if $i\leq -1$ and $\Fil^iK(1)_0=0$ if $i>-1$.

If $D$ is a filtered $(\Phi,N)$-module, we denote $D(1)=D\otimes K(1)$. The trivialisation $\psi$ of the Tate twist permits to describe an isomorphic module $D(1)_0$ having the same underlying $K_0$-vector-space as $D$, the same $N$, the Frobenius $\Phi'=(1/p)\Phi$ and the filtration translated by $-1$:  $\Fil^iD(1)_0=\Fil^{i+1}D$.
\end{eexm}

\subsection{Transverse monodromy}
In all this paragraph, we use the trivialisation $\psi\colon \qp(1)\to \qp(1)_0$ of the Tate twist and the induced isomorphism $i_\psi\colon K(1)\to K(1)_0$ as in \ref{exmtate}.

For certain particular $p$-adic representations of a $p$-adic field $K$, and notably for those associated to semi-stable abelian varieties over $K$, the monodromy is visible over $\qp$ directly on the representation. This is due to the fact that, in those cases, the monodromy gives a morphism in the category of filtered $(\Phi,N)$ modules, which implies that it is visible in any realisation of this category, and notably in the etale $p$-adic realisation. Let us define those representations.

\begin{edefin}
We say that an object $D$ of $\MF_K^a(\Phi,N)$ has \textrm{transverse monodromy} if
$$N\Fil^iD_K\subset \Fil^{i-1}D_K.$$
We denote by $\MF_K^{a,tm}(\Phi,N)$ the subcategory of such objects in $\MF_K^a(\Phi,N)$.
\end{edefin}

We will show that this condition is verified for the filtered module of an
abelian variety. A counter example to the general case is given in the
mathscinet review of the article \cite{Breuil}, in the case of the Galois
representation associated to a modular form.

\begin{eprop}
$\MF_K^{a,tm}(\Phi,N)$ is a tannakian subcategory of $\MF_K^a(\Phi,N)$.
\end{eprop}
\begin{proof}
It suffices to verify that transversality of monodromy is stable by tensor product, subobject, quotient and dual in the category of admissible filtered $(\Phi,N)$-modules.
The verification for tensor products and dual is straight forward.
Let $D$ such a module with transverse monodromy and $D'\subset D$ be a subobject of $D$. Since morphisms are strict (this is a consequence given in \cite{Fontaine}, 4.3.3 of the implication "admissible$\Rightarrow$weakly-admissible" \cite{Fontaine}, 5.4.2), we have the relation $\Fil^iD'_K=D'_K\cap \Fil^iD_K$. This permits us to obtain
$$N\Fil^iD'_K=N(D'_K\cap \Fil^iD_K)\subset D'_K\cap \Fil^{i-1}D_K=\Fil^{i-1}D'_K$$
which is the required result for submodules. To finish, the transversality for quotients is a consequence of the transversality for duals and submodules.
\end{proof}

We denote $\Rep_{st,tm}(\GKK)$ the essential image of $\MF_K^{a,tm}(\Phi,N)$ by  $\V_{st}$ and $\Rep_{pst,tm}(\GKK)$ the category of potentially semi-stable representations such that on every finite extension $L$ of $K$ on which the Galois representation is semi-stable, it belongs to $\Rep_{st,tm}(\GLL)$, where $\GLL:=\Gal(\overline{L}/L)$. We remark that since $\MF_K^{a,tm}(\Phi,N)$ is a tannakian sub-category of $\MF_K^a(\Phi,N)$, it follows that $\Rep_{st,tm}(\GKK)$ (resp. $\Rep_{pst,tm}(\GKK)$) is a tannakian sub-category of $\Rep_{st}(\GKK)$ (resp. $\Rep_{pst}(\GKK)$).

Denote $\MF_K^a(\Phi,N)(N_{st})$ the category of couples of an admissible filtered $(\Phi,N)$-module $D$ and a nilpotent morphism $N_{st}\colon D(1)\to D$. By nilpotent, we mean that for some integer $n$, the composition
$$
(N_{st}\otimes \id_{K(n-1)})\circ
\cdots\circ(N_{st}\otimes \id_{K(i)})\circ
\cdots\circ N_{st}\colon D(n)\to D,
$$
is the zero morphism, where $K(i):=K(1)^{\otimes i}$ for $i\in\N$.

\begin{elem}
\label{functor-monod-st}
Let $D$ be an admissible filtered $(\Phi,N)$-module. The monodromy map
$N:D\to D$ induces a $K_0$-linear application $N_{st,\psi}\colon D(1)\to D$
that depends on the trivialization $\psi$ of the Tate twist. Moreover,
the map $N_{st,\psi}$ is a morphism of filtered $(\Phi,N)$-module if and only if $D$ has transverse monodromy. This gives a tensor functor
$$m_{st,\psi}\colon \MF_K^{a,tm}(\Phi,N)\to \MF_K^a(\Phi,N)(N_{st}).$$
\end{elem}
\begin{proof}
Let $D$ be an admissible filtered $(\Phi,N)$-module.
We still denote $N$ the linear map $D(1)_0\to D$ given by $N$. We now show that it is a morphism of filtered $(\Phi,N)$-module. This application commutes with the action of $\Phi$ since $N.\Phi'=1/p.p.\Phi.N=\Phi.N$. It commutes with $N$ since $N^2=N^2$. Those two properties are showed without hypothesis on $D$. The transversality property now appears. $N$ is a morphism of filtered $(\Phi,N)$-module if and only if
$$\forall i, N\Fil^iD(1)_{0,K}\subset \Fil^iD_K,$$
which is equivalent to
$$\forall i, N\Fil^iD_K\subset \Fil^{i-1}D_K,$$
i.e. to the fact that $D$ has transverse monodromy.
We thus have showed that $N$ induces a morphism $N\colon D(1)_0\to D$ in $\MF_K^a(\Phi,N)$ if and only if $D$ has transverse monodromy. This monodromy $N\colon D(1)_0\to D$ is functorial in $D$ and compatible with the tensor product.
Using the isomorphism $\id_D\otimes i_\psi\colon D(1)\to D(1)_0$, we obtain a morphism $N_{st,\psi}=N\circ (\id_D\otimes i_\psi)\colon  D(1)\to D$ that depends on $t_\psi=\psi^{-1}(1)$ up to a scalar. This morphism $N_{st,\psi}$ is called the \emph{monodromy morphism} of $D$. The rule $D\mapsto (D,N_{st,\psi})$ is functorial.

We now show that the monodromy morphism $N_{st,\psi}\colon D(1)\to D$ is compatible with the tensor product. If we choose two filtered $(\Phi,N)$-modules with transverse monodromy $D_1$ and $D_2$, then the monodromy morphism of the module $D=D_1\otimes D_2$ is given by 
$$
N_{st,\psi}=\id_{D_1}\otimes N_{2,st,\psi}+N_{1,st,\psi}\otimes \id_{D_2}
$$
where $N_{1,st,\psi}\colon D_1(1)\to D_1$ (resp. $N_{2,st,\psi}\colon D_2(1)\to D_2$, resp. $N_{st,\psi}\colon (D_1\otimes D_2)(1)\to (D_1\otimes D_2)$) is the monodromy morphism for $D_1$ (resp. $D_2$, resp. $D=D_1\otimes D_2$).
So let $N\colon D(1)_0\to D$, $N_1\colon D_1(1)_0\to D_1$ and $N_2\colon D_2(1)_0\to D_2$ denote the respective monodromy morphisms with trivialized Tate twist. We develop
$$
\begin{array}{cl}
N_{st,\psi} & =N\circ (\id_D\otimes i_\psi)\\
         & =(N_1\otimes \id_{D_2}+\id_{D_1}\otimes N_2)\circ ((\id_{D_1}\otimes \id_{D_2})\otimes i_\psi)\\
	 & =(N_1\circ (\id_{D_1}\otimes i_\psi))\otimes \id_{D_2}+\id_{D_1}\otimes (N_2\circ (\id_{D_2}\otimes i_\psi))
\end{array}
$$
to obtain the desired equality.
We thus have defined a tensor functor
$$
\begin{array}{cccc}
m_{st,\psi}\colon  & \MF_K^{a,tm}(\Phi,N) & \to     & \MF_K^a(\Phi,N)(N_{st})\\
             & D                    & \mapsto & (D,N_{st,\psi}\colon D(1)\to D)
\end{array}.
$$
\end{proof}

Since the category $\MF_K^{a,tm}(\Phi,N)$ is equivalent to the category $\Rep_{st,tm}(\GKK)$, for every $p$-adic semi-stable representation $V$ of $\GKK$ the associated module of which has transverse monodromy, $N_{st,\psi}$ induces a morphism
$N_{p,\psi}\colon V(1)\to V$ of $p$-adic representations.

Denote $\Rep_{st}(\GKK)(N_p)$ the category whose objects are couples $(V,N_p)$ of a semi-stable representation and a nilpotent morphism $N_p\colon V(1)\to V$, and whose morphisms are the morphisms of representations that commute with the given $N_p$.

\begin{edefin}
\label{functor-padic-monod}
We have defined a functor
$$
\begin{array}{cccc}
m_{p,\psi}\colon  & \Rep_{st,tm}(\GKK) & \to     & \Rep_{st}(\GKK)(N_p)\\
            & V                  & \mapsto & (V,N_{p,\psi})
\end{array}.
$$
\end{edefin}

This functor is in fact a tensor functor because it is given by the following diagram
$$
\xymatrix{
\Rep_{st,tm}(\GKK)\ar[r]^{\D_{st}}\ar[d]^{m_{p,\psi}} & \MF_K^{a,tm}(\Phi,N)\ar[d]^{m_{st,\psi}}\\
\Rep_{st}(\GKK)(N_p)                      & \MF_K^a(\Phi,N)(N_{st}) \ar[l]^{\V_{st}}}
$$
that is commutative up to an isomorphism of functors.

\subsection{Exponential of the $p$-adic monodromy}
\label{expo-monod-padic}
We will now use the usual formalism of Tannakian categories given in \cite{De1}, II. A useful example in this paragraph is the following.
\begin{eexm}
\label{exm-tannakian}
The category $\Vect_\qp(N)$ of couples $(V,N)$ of a $\qp$-vector space $V$ and a $\qp$-linear nilpotent endomorphism $N\colon V\to V$ endowed with the tensor product
$$(V_1,N_1)\otimes (V_2,N_2)=(V_1\otimes V_2,N_1\otimes \id_{V_2}+\id_{V_1}\otimes N_2)$$
is a tannakian category for which the forgetful functor $o\colon (V,N)\mapsto V$ is a fibre functor.
The associated group scheme $\Aut^\otimes(o)$ is the additive group $\G_{a,\qp}$. This gives a tensor equivalence of categories
$$o\colon \Vect_\qp(N)\to\Rep(\G_{a,\qp}).$$
\end{eexm}

Recall that we have chosen a trivialisation $\psi\colon \qp(1)\to\qp(1)_0$ of the Tate twist. If $V$ is a $p$-adic representation in $\Rep_{st,tm}(\GKK)$, we have defined a morphism $N_{p,\psi}\colon V(1)\to V$ of $p$-adic representations. Now define $V(1)_0$ to be $V\otimes \qp(1)_0$. The underlying vector space of $V(1)_0$ is canonically isomorphic to $V$. We obtain a morphism $N_\psi\colon V(1)_0\to V$ by setting $N_\psi:=N_{p,\psi}\circ (\id_{V}\otimes \psi^{-1})$. This morphism gives an endomorphism of the underlying $\qp$-vector space of $V$.

\begin{eprop}
\label{expsttm}
The rule
$$
\begin{array}{cccc}
\omega_\psi\colon  & \Rep_{st,tm}(\GKK) & \to     & \Vect_\qp(N)\\
	     & V                  & \mapsto & (V,N_\psi=N_{p,\psi}\circ (\id_{V}\otimes \psi^{-1}))
\end{array}
$$
gives a tensor functor to the category of vector spaces endowed with a nilpotent endomorphism and we consequently have a morphism
$$e^N_\psi\colon \G_{a,\qp}\to \Aut^\otimes(\omega_p)$$
from the additive group to the tannakian fundamental group $\Aut^\otimes(\omega_p)$ associated to the natural fibre functor $\omega_p=o\circ \omega_\psi$ on the category $\Rep_{st,tm}(\GKK)$. This morphism is called an \emph{exponential of the $p$-adic monodromy up to the Tate twist}.
\end{eprop}
\begin{proof}
We must first verify that if $f\colon V_1\to V_2$ is a morphism in the category $\Rep_{st,tm}(\GKK)$, the associated linear map over $\qp$ commutes with $N_\psi$. So we must show that the outer rectangle of following diagram commutes.
$$
\xymatrix{
V_1\ar[d]^{f}\ar[r]^{\id_{V_1}\otimes\psi^{-1}} & V_1(1) \ar[d]^{f(1)}\ar[r]^{N_{1,p,\psi^{-1}}} & V_1\ar[d]^{f}\\
V_2          \ar[r]^{\id_{V_2}\otimes\psi^{-1}} & V_2(1)              \ar[r]^{N_{2,p,\psi}} & V_2}
$$
The right square of the diagram is commutative because the map $V\mapsto (V,N_p)$ is a functor $m_{p,\psi}\colon \Rep_{st,tm}(\GKK)\to\Rep_{st}(\GKK)(N_p)$, as remarked in \ref{functor-padic-monod}. The left square of the diagram is obviously commutative so we have proved that the whole diagram is commutative, which shows that we really have a functor.

To verify that it is a tensor functor, we must verify that the trivialisation of the Tate twist does not affect the tensor product. So let $V_1$ and $V_2$ be two objects of $\Rep_{st,tm}(\GKK)$ (so their image by $m_{p,\psi}$ is endowed with $N_{1,p,\psi}\colon V_1(1)\to V_1$ and $N_{2,p,\psi}\colon V_2(1)\to V_2$). Then since the functor $m_{p,\psi}\colon V\mapsto (V,N_{p,\psi})$ is a tensor functor (see \ref{functor-padic-monod}), $W=V_1\otimes V_2$ is a Galois representation endowed with a nilpotent morphism $N_{p,\psi}\colon W(1)\to W$ that is given by $N_{p,\psi}=N_{1,p,\psi}\otimes \id_{V_2}+\id_{V_1}\otimes N_{2,p,\psi}$. We want to show that trivialisation of the Tate twist is compatible with the tensor product of monodromy. So we must show that
$$N_{\psi}=N_{1,\psi}\otimes \id_{V_2}+\id_{V_1}\otimes N_{2,\psi}.$$
Developing the left term
$$
\begin{array}{cl}
N_{\psi} & =N_{p,\psi}\circ (\id_W\otimes \psi)\\
         & =(N_{1,p,\psi}\otimes \id_{V_2}+\id_{V_1}\otimes N_{2,p,\psi})\circ ((\id_{V_1}\otimes \id_{V_2})\otimes \psi)\\
	 & =(N_{1,p,\psi}\circ (\id_{V_1}\otimes \psi))\otimes \id_{V_2}+\id_{V_1}\otimes (N_{2,p,\psi}\circ (\id_{V_2}\otimes \psi))
\end{array}
$$
gives exactly the required equality.
So we have a tensor functor, since
$$\omega_\psi(V_1\otimes V_2)=\omega_\psi(V_1)\otimes \omega_\psi(V_2).$$
The rest of the proposition follows from \ref{exm-tannakian} and the usual
tannakian formalism (see \cite{De1}, II).
\end{proof}

\begin{ecor}
\label{expVsttm}
Let $V$ be an object of $\Rep_{st,tm}(\GKK)$ and let $H_{V,p}$ denote the Zariski closure of the image of this Galois representation. We have a natural morphism
$$e^N_\psi\colon \G_{a,\qp}\to H_{V,p}$$
that gives an ``exponential of the $p$-adic monodromy up to a Tate twist''.
\end{ecor}
\begin{proof}
Let $\omega_p$ be the natural fibre functor on $\Rep_{st,tm}(\GKK)$, i.e. the functor that associates to any representation its underlying $\qp$-vector space.
Let $V$ be an object of $\Rep_{st,tm}(\GKK)$ and $\langle V\rangle$ be the tannakian subcategory of $\Rep_{st,tm}(\GKK)$ generated by $V$. We denote by $\omega_p|_{\langle V\rangle}$ the restriction of $\omega_p$ to this subcategory. By the usual tannakian machinery (see for example \cite{Fontaine}, 1.2.3), we have an isomorphism
$$\Aut^\otimes(\omega_p|_{\langle V\rangle})\cong H_{V,p}$$
and a natural projection associated to the tensor functor $\langle V\rangle\to \Rep_{st,tm}(\GKK)$
$$\Aut^\otimes(\omega_p)\twoheadrightarrow \Aut^\otimes(\omega_p|_{\langle V\rangle}).$$
This permits to define
$$e^N_\psi\colon \G_{a,\qp}\to \Aut^\otimes(\omega_p)\to H_{V,p}$$
by composition.
\end{proof}

\subsection{$p$-adic monodromy for $\mathrm{H}_1$ and good reduction}
If a $p$-adic representation $V$ is the first etale $p$-adic cohomology group $\mathrm{H}^1$ of a semi-stable projective variety over $K$, the associated filtered $(\Phi,N)$-module has a filtration with only nontrivial terms
$$
\begin{array}{cl}
\Fil^0\D_{st,K}(V) & = \D_{st,K}(V),\\
\Fil^1\D_{st,K}(V) &\subset \D_{st,K}(V),\\
\Fil^2\D_{st,K}(V) & =  0.
\end{array}
$$
For such a module $D$, the condition of transversality of monodromy
$$N\Fil^iD_K\subset \Fil^{i-1}D_K$$
is trivially verified, so there is an exponential of the monodromy morphism as in \ref{expVsttm}. This remark applies to abelian varieties. This gives a $p$-adic counterpart, in the case of abelian varieties, to the construction of Raynaud \cite{Ray2}, p315 of the $\ell$-adic monodromy for $1$-motives for $\ell\neq p$.

\begin{eprop}[$p$-adic criterion for good reduction]
\label{padic-criteria}
Let $A$ be an abelian variety over $K$ and $\rho\colon \GKK\to \GL(V_p(A))$ be the corresponding $p$-adic representation. If the Zariski closure $H_{V,p}$ of the image of $\rho$ does not contain any unipotent $\qp$-rational point of index $2$, then $\rho$ is potentially cristalline and $A$ has potentially good reduction.
\end{eprop}
\begin{proof}
We can suppose that $\rho$ is semi-stable. We denote $D(A)$ the filtered $(\Phi,N)$-module $\D_{st,K}(V_p(A))$.
The monodromy $N_{st,\psi}\colon D(A)(1)\to~D(A)$ in the category $\MF_{K}^{a,tm}(\Phi,N)$ is nilpotent of index $2$, i.e. $N_{st,\psi}\circ N_{st,\psi}(1)=0$. This comes from the fact that if $N\colon D(A)\to D(A)$ denotes the monodromy morphism, $N^2=0$. We will prove this using the $1$-motive $M/K$ constructed by Raynaud in \cite{Ray2}, 4.2.2. This is a complex $M:=[X\to G]$ with $X$ (in degree $-1$) a twisted free $\Z$-module and $G$ (in degree $0$) an extension of an abelian variety $B$ by a torus $T$. This $1$-motive has the same de Rham and $p$-adic realisation as $A$.
It has a weight filtration given by $W_{-3}=[0\to 0]$, $W_{-2}=[0\to T]$, $W_{-1}=[0\to G]$ and $W_0=[X\to G]$ with successive quotients $\Gr_0M=X[-1]$, $\Gr_{-1}M=B$ and $\Gr_{-2}M=T$.
We denote $M_{dR}$ the de Rham realization $H^1_{dR}(M)$. We see in \cite{CoIo}, 2.1 that there is a natural map
$\Gr_0M_{dR}\to \Gr_{-2}M_{dR}$ and Coleman proves in \cite{Co1} that it is induced by Grothendieck's monodromy pairing given in \cite{SGA7}, IX, Th 10.4.
The monodromy $N\colon M_{dR}\to M_{dR}$ is given by the sequence
$$N\colon M_{dR}\to\Gr_0M_{dR}\to\Gr_{-2}M_{dR}=W_{-2}M_{dR}\to M_{dR}.$$
We now remark that the composition
$$\Gr_{-2}M=W_{-2}M\to M\to \Gr_{0}M$$
is clearly trivial and so is its de Rham realization.
This proves that $N\circ N=0$, so $N_{st,\psi}\circ N_{st,\psi}(1)=0$.

Since we have an equivalence of additive categories, this nilpotency of index $2$ is also verified for the $p$-adic monodromy morphism
\mbox{$N_{p,\psi}\colon V_p(A)(1)\to V_p(A)$}. Using the trivialisation of the Tate twist $\qp(1)\overset{\psi}{\to}\qp(1)_0$, we obtain a $\qp$-linear map $N_{\psi}\colon V_p(A)\to V_p(A)$, as in \ref{expsttm} that is nilpotent of index $2$, i.e. $N_{\psi}^2=0$. The exponential of the monodromy morphism $e^N_{\psi}\colon \G_{a,\qp}\to \GL(V_p(A))$ is given by the image of the universal point $T=\id_{\G_a}\in \G_a(\qp[T])$ of the additive group: 
$$
\begin{array}{cccc}
e^N_{\psi}\colon  & \G_a(\qp[T]) & \to     & \GL(V_p(A))(\qp[T])\\
            & T          & \mapsto & \sum_{i\geq 0}\frac{N_{\psi}^iT^i}{i!}
\end{array}.
$$
Since $N_{\psi}$ is nilpotent of index $2$, the image of this morphism is unipotent of index $2$, i.e. every element $g$ of this image verifies the equality $(g-\id)^2=1$, calculated in the $\qp$-linear endomorphisms of $V$. By \ref{expVsttm}, this morphism factors through the Zariski closure $H_{V,p}$ of the image of $\rho$. The image of $e^N_{\psi}$ in $H_{V,p}$ is trivial because $H_{V,p}(\qp)$ contains no nontrivial unipotent element of index $2$ in $\GL(V_p)$. This implies that $N_{\psi}$ and therefore $N_{p,\psi}\colon V_p(A)(1)\to V_p(A)$ and $N_{st,\psi}\colon D(A)(1)\to D(A)$ are trivial. So the original monodromy $N$ of the filtered $(\Phi,N)$-module $D(A)$ is trivial and the representation $V_p(A)$ is potentially cristalline.
The work of Coleman-Iovita \cite{CoIo}, Th1, p173 or Breuil \cite{Breuil2}, Cor 1.6, p493 shows that the Galois representation corresponding to an abelian variety is potentially cristalline if and only if the abelian variety has potentially good reduction.
\end{proof}

\subsection{The criterion}
We cautiously inform the reader that $K$ will now denote a number field.
\begin{etheo}
\label{CRITERION}
Let $K\hookrightarrow \C$ be a number field embedded in $\C$ and $A$ be an abelian variety over $K$. Let $V=\Hun(A_{\C}^{an},\Q)$ be its first Betti homology group and $\rho\colon \MT(A)\to \GL(V)$ be the natural representation of its Mumford-Tate group. Suppose there exists a prime $\ell$ such that the image of
$$\rho_{(\ql)}\colon \MT(A)(\ql)\to \GL(V)(\ql)$$
does not contain any unipotent of index $2$,
then $A$ has potentially good reduction at any finite place of $K$. 
\end{etheo}
\begin{proof}
By the comparison theorem between etale and classical cohomology (see \cite{SGA4}, XI) and the comparison between the dual of the first etale cohomology group and the Tate module (see Milne \cite{Mi2}, 15.1), we know that there is a canonical isomorphism $V\otimes_\Q \ql\cong T_\ell(A)\otimes_{\Z_\ell}\Q_\ell$.
It results from theorems of Deligne about absolute Hodge cycles \cite{De1}, Exp I, 2.9,2.11 that there exists a finite extension $K'$ of $K$ such that the Galois representation on the $\ell$-adic Tate module of $A$ factors through the group of $\ql$-values of the Mumford-Tate group, i.e. there is a morphism
$$\rho'\colon \Gal(\overline{K}/K')\to \MT(A)(\ql)\overset{\rho_{(\ql)}}{\to} \GL(V)(\ql)\cong \GL(T_\ell(A)\otimes_{\Z_\ell}\Q_\ell)$$
that gives a factorisation of the usual Galois representations on $\ell^n$ torsion points of $A$ for all $n>0$ through $\rho_{(\ql)}$.

Let $v$ be a finite place of $K$.
Suppose first that $\ell$ doesn't divide $v$. Then, by the monodromy theorem of Grothendieck \cite{SGA7}, Exp I, 3.6, there is a finite extension $L$ of $K_v$, an extension $w$ of $v$ to $L$ and a place $\bar{w}$ of $\overline{L}$ over $w$ such that the image of the inertia subgroup $\II(\bar{w})$ of $\DD(\bar{w})$ is unipotent of index $2$ in $\GL(T_\ell(A)\otimes_{\Z_\ell}\Q_\ell)$. More precisely, we have $(\rho'(\sigma)-\id)^2=0$ in $\End(V_{\ql})$ for all $\sigma$ in $\II(\bar{w})$. Since the image of the representation $\rho_{(\ql)}$ of the Mumford-Tate group over $\ql$ contains no nontrivial unipotent of index $2$, the image of the inertia group $\II(\bar{w})$ must be trivial. So the original Galois representation was potentially unramified at $v$. The N\'eron-Ogg-Shafarevich criterion \cite{SeTa}, Th 2, gives the potentially good reduction of $A$ at $v$.

Suppose now that $\ell$ divides $v$. Then we apply our $p$-adic criterion \ref{padic-criteria} to the representation of the decomposition group $\DD(\bar{v})\subset \GL(T_\ell(A)\otimes_{\Z_\ell}\Q_\ell)$, which concludes the proof.
\end{proof}

\section{Classification of representations satisfying \ref{CRITERION}}
\label{combinpart}
Let $K\hookrightarrow \C$ be a number field and $A$ be an abelian variety over $K$. Let $\rho\colon G\to \GL(V)$ be the natural representation of its Mumford-Tate group. Let $\ell$ be a prime and $\rho_\ql\colon G_\ql\to \GL(V)_\ql$ be the extension of scalars of $\rho$ to $\ql$.
We will first explain in simple terms what kind of representations $\rho$ and $\rho_\ql$ arise this way. This gives a combinatorial datum called a \emph{polymer}.

After that, we will use this datum to give a criterion for the representation \mbox{$\rho\colon G\to\GL(V)$} to have the following property:
for some prime $\ell$, the image of the representation
$$\rho_{(\ql)}\colon G(\ql)\to \GL(V)(\ql)$$
of $\ql$-points of $G$ induced by $\rho$ contains no nontrivial unipotents of index $2$.

\subsection{Classification of representations of Mumford-Tate groups}
\subsubsection*{Unipotents in reductive groups}
\label{unip-reductive}
We will not go into the details of the classification of representations of reductive groups here.
We refer to \cite{SGA3III}, \cite{Jantzen} and \cite{Spr} for the structure theory of reductive groups. Since our work focuses on reductive groups over non algebraically closed fields, we use the language of group schemes.
Let $F$ be a field of characteristic $0$, $G$ be a reductive group over $F$ and $T\subset G$ be a maximal torus defined over $F$.

To these data are associated two twisted free abelian group schemes of finite type over $F$ (or equivalently two free abelian groups of finite type equipped with continuous-actions of $\Gal(\overline{F}/F)$)
$$X^*(T)=\underline{\Hom}(T,\G_m)\textrm{ and }X_*(T)=\underline{\Hom}(\G_m,T)$$
and two closed sub-schemes
$$R=R(G,T)\subset X^*(T)\textrm{ and }R^\vee=R^\vee(G,T)\subset X_*(T)$$
called the \emph{scheme of roots (resp. of coroots)} of $(G,T)$.
There is a perfect pairing
$$\langle\cdot,\cdot\rangle\colon X^*(T)\times X_*(T)\to \Z$$
given by composition of characters with cocharacters and by the isomorphism $\rm{End}(\G_m)\cong \Z$.
The scheme of roots $R$ is the scheme of weights of the maximal torus in the adjoint representation (see \cite{SGA3III}, Exp XIX, 3.8). For $g\in G$, we denote $\ad_g\colon \Lie(G)\to\Lie(G)$ the tangent map corresponding to $h\mapsto ghg^{-1}$.
Let
$$
\begin{array}{cccccc}
\rho_{ad}\colon & T & \hookrightarrow & G  & \to     & \GL(\mathrm{Lie}(G))\\
           	&   &                 & g  & \mapsto & \ad_g
\end{array}
$$
denote the adjoint representation. If $r\in X^*(T)(F)$ and $A$ is an $F$-algebra, define the \emph{weight scheme of $r$ in $\Lie(G)$} by
$$\Lie(G)^r(A):=\left\{
\begin{array}{c|l}
h\in \Lie(G)\otimes_F A &\rho_{ad,A'}(t)\cdot h=r_{A'}(t)\cdot h,\\
			&\forall t\in T(A'),A\to A'
\end{array}\right\}$$
Let now
$$
R(A):=\left\{r\in X^*(T_A)\,|\,\forall s\in\Spec(A), r_s\neq 0\textrm{ and } \Lie(G_A)^r(s)\neq 0\\
\right\}.$$
Since it isn't crucial for our work, we refer to \cite{SGA3III}, Exp XXII for the definition of the scheme of coroots in the non split case and to \cite{Jantzen}, II, 1.3 in the split case.
\begin{edefin}[Root datum]
The data $\Psi=(X^*(T),R,X_*(T),R^\vee)$ of those four schemes and of the pairing $\langle\cdot,\cdot\rangle$ is called \emph{the root datum of the couple $(G,T)$}.
\end{edefin}
The root datum plays a crucial role in the classification of reductive groups: it allows to classify split reductive groups (i.e. groups with $T\cong \G_m^r$) up to isomorphism.
It also gives precise information about unipotent elements in non-split reductive groups. Indeed, for each root $r\in~R(F)$, there exists a (unique up to an element of $F^*$) morphism $u_r\colon \G_a\to G$ (see \cite{SGA3III}, Chap XXII, 1.1) satisfying for all $F$-algebra $A$
$$tu_{r,A}(x)t^{-1}=u_{r,A}(r_{A}(t)\cdot x),\forall t\in T(A),x\in A.$$
The image of this morphism is called the \emph{standard unipotent subgroup corresponding to the root $r$}.

We also recall some basic properties of unipotent elements in reductive groups from a functorial point of view.
For now, denote $G$ a linear algebraic group over a field $F$ of characteristic $0$.
Define the functor $\underline{\Hom}(\G_a,G)$ on $F$-algebras by $$\underline{\Hom}(\G_a,G)(A)=\Hom_{A-Gr}(\G_{a,A},G_A)$$
for every $F$-algebra $A$.
There is a natural injective morphism
$$
\begin{array}{ccc}
\underline{\Hom}(\G_a,G) & \to     & G\\
(u\colon \G_a\to G)      & \mapsto & u(1)
\end{array}
$$
and its image $\U(G)$ is \emph{the subscheme of unipotent elements in $G$}
(every unipotent element $U$ in $G$ induces a morphism $u:\G_a\to G$  such that $u(1)=U$ by $u(T)=\exp(T\log(U))$ for $T=\id\in \G_a(F[T])$).

The fact that this is a closed subscheme can be verified on a general linear group $\GL(V)$ since every linear algebraic group admits a closed embedding in such a group. In the case of $\GL(V)$ we know that $\U(\GL(V))$ is given by $u\in\GL(V)$ satisfying $(u-\id)^{dim(V)}=0$. This is a polynomial condition and gives a closed subscheme.

We now recall a well known fact.
\begin{eprop}
\label{prop-unip-adjoint}
If $G$ is reductive, and $\tilde{G}$ denotes the universal covering of $G^{ad}$, the sequence
$\tilde{G}\to G\to G^{ad}$ gives isomorphisms
$$\U(\tilde{G})\cong \U(G)\cong \U(G^{ad})$$
on unipotent subschemes.
\end{eprop}
\begin{proof}
The exponential map gives an isomorphism between $\U(G)$ and nilpotent elements
in $\Lie(G)$ and the result is trivial on the level of Lie algebras.
\end{proof}

\subsubsection*{Admissible minuscule representations}
We now recall the definition of minuscule representations that best suits our work. We refer to Boubaki \cite{Boulie78}, Chap VIII, \S7, n$^\circ$ 3, for the treatment of this notion in the context of Lie algebras. There is a little mistake\footnote{For the $B_l$ type, replace $\omega_1$ by $\omega_l$ and for the $C_l$ type, replace $\omega_l$ by $\omega_1$} in the list of minuscule weights in this reference, so we use \cite{Se7}, Annexe, as a reference for the list of minuscule weights.
\begin{edefin}[Minuscule representation]
\label{definminuscule}
Suppose $F$ is algebraically closed and $G$ is semi-simple. An irreducible representation $V$ of $G$ is called \emph{minuscule} if for every root $r\in R(F)$, the image of the standard unipotent subgroup $u_r\colon \G_a\to G$ in $\GL(V)$ is unipotent of index $2$.
\end{edefin}

We still suppose $F$ algebraically closed.
Let $G/F$ be a reductive group endowed with a maximal torus $T$ and a Borel subgroup $B$. These data give a root datum $\Psi$ with a set $R^+\subset R$ of positive roots and a basis $S\subset R^+$ of simple roots. Simple roots correspond bijectively with  vertices of the Dynkin diagram $D$ of $G$.

We will say that $G$ is \emph{semi-simple simply connected} if $X^*(T)$ admits a basis $(\omega_\alpha)_{\alpha\in S}$ dual to the basis $(\alpha^\vee)_{\alpha\in S}$ of $R^\vee$. The $\omega_\alpha$ for $\alpha\in S$ are called fundamental weights. They also correspond bijectively with vertices of the Dynkin diagram.
We now recall the basics of representation theory of reductive groups, that are contained in \cite{Jantzen}, Part II.
Set
$$X^*(T)_+=\{\lambda\in X^*(T)\,|\,\langle\lambda,\alpha^\vee\rangle\geq 0\textrm{ for all } \alpha\in S\}.$$
The elements of $X^*(T)_+$ are called the dominant weights of $T$ with respect to $R^+$.
\begin{etheo}
\label{domheighest}
Keeping the notations above, and if we suppose that $G$ is semi-simple simply connected, there is a bijection between the set $X^*(T)_+$ of dominant weights of $T$ with respect to $R^+$ and the set of isomorphism classes of finite dimensional simple $G$-modules. A dominant weight corresponding to some module $V$ by this bijection will be called \emph{the highest weight of $V$}.
\end{etheo}

The fundamental weights $(\omega_\alpha)_{\alpha\in S}$ are dominant.
This implies, using theorem \ref{domheighest}, that if $G$ is semi-simple simply connected, every fundamental weight is indeed the highest weight of a simple $G$-module.

The weights corresponding to minuscule representations of the simply connected covering of the adjoint group of $G$ are called \emph{minuscule weights}. We know from Bourbaki (cf. the end of \cite{Boulie78}, Chap VIII, \S7, n$^\circ$ 3) that minuscule weights are fundamental. Thus they correspond to some vertices of the Dynkin diagram of $G$.

Depending on the type of $G$ in the classification of reductive groups, the set of minuscule weights is given in \cite{Se7}, Annexe. We now recall the list of weights that are admissible for the representations that will be the building blocks of Hodge type representations. They are given in the table \cite{De4}, 1.3.9.
\begin{edefin}[Admissible minuscule weight]
\label{adminuscule}
Let $\Psi=(X^*(T),R,X_*(T),R^\vee)$  be the root datum of a simple semi-simple simply connected group $(G,T)$, for which we fix a system of positive roots $R^+\subset R$ and a basis of simple roots $S=\{\alpha_1,\dots,\alpha_n\}\subset R^+$. Let $(\omega_i)_{i=1,\dots,n}$ be the corresponding basis of fundamental weights. Then the set of \emph{admissible minuscule weights for $\Psi$ and $R^+$} is given by:
\begin{itemize}
\item for $\Psi$ of type $A_n$ ($n\geq 1$), $\{\omega_1,\omega_n\}$,
\item for $\Psi$ of type $B_n$ ($n\geq 2$), $\{\omega_n\}$,
\item for $\Psi$ of type $C_n$ ($n\geq 2$), $\{\omega_1\}$,
\item for $\Psi$ of type $D_n$ ($n\geq 4$), $\{\omega_1,\omega_{n-1},\omega_n\}$,
\item for $\Psi$ of other types, $\emptyset$.
\end{itemize}
The corresponding vertices of the Dynkin diagram are called \emph{admissible minuscule vertices}.
\end{edefin}
We see from \cite{Se7}, Annexe, that admissible minuscule weights are indeed minuscule.

\subsubsection*{Polymers and Hodge type representations}

Let's come back to a reductive group $G$ over $F$ with $F$ not necessarily algebraically closed. We denote $G^{ad}$ the adjoint group of $G$ and $\tilde{G}$ the universal covering of $G^{ad}$. We have a sequence $\tilde{G}\to G\to G^{ad}$. We decompose $\tilde{G}$
$$\tilde{G}_{\overline{F}}\cong \prod_{i\in I}G_i$$
as a product of simple groups over $\overline{F}$.
The set $I$ of simple factors of $\tilde{G}_{\overline{F}}$ is endowed with a $\Gal(\overline{F}/F)$-action. Let $D$ be the Dynkin diagram of $\tilde{G}_{\overline{F}}$, endowed with its natural $\Gal(\overline{F}/F)$-action. There is a natural equivariant projection $\pi\colon D\to I$. We denote $\mathscr{P}(D)$ the set of parts of $D$ and $D=\coprod_{i\in I} D_i$ the decomposition of the Dynkin diagram in connected components. As explained above, every vertex of the Dynkin diagram $D_i$ corresponds to a fundamental weight, i.e. to the heighest weight of a fundamental representation of the factor $G_i$ (that is semi-simple simply connected).

We are now able to define the kind of representations that appear in our situation and the corresponding combinatorial datum, called a \emph{polymer} (after Addington \cite{Add}).
\begin{edefin}[Polymer]
\label{defin-polymer}
Let $G$ be a reductive group over $F$. A \emph{polymer for $G$} is a set
$\mathscr{S}\subset~\mathscr{P}(D)$ of subsets of the Dynkin diagram of $\tilde{G}$ that verifies:
\begin{enumerate}
\item $\mathscr{S}$ is $\Gal(\overline{F}/F)$-stable,
\item $\mathscr{S}$ is covering on the set $I$ of simple factors of $\tilde{G}$, i.e. $\underset{T\in\mathscr{S}}{\cup}\pi(T)=I$,
\item if $T\in~\mathscr{S}$ and $T\cap D_i\neq \emptyset$ then $T\cap D_i$ is reduced to an admissible minuscule vertex of the connected component $D_i$ of $D$ (in the sense of \ref{adminuscule}).
\end{enumerate}
\end{edefin}

\begin{edefin}[Hodge type representation]
\label{rephodgetype}
A \emph{Hodge type representation} of $G$ is a representation
$$G\to \GL(V)$$
such that there exists a polymer $\mathscr{S}(V)$ for $G$ such that every irreducible component $W$ of the restriction of $V_{\overline{F}}$ to $\tilde{G}_{\overline{F}}$ is isomorphic to a tensor product
$\otimes_{s\in T}W_s$
for some $T\in\mathscr{S}(V)$. In this tensor product, $W_s$ denotes the minuscule representation of the simple factor $G_{\pi(s)}$ of $\tilde{G}_{\overline{F}}$ of heighest weight $\omega_s$ corresponding to the unique $s\in T\cap D_{\pi(s)}$ given by definition of a polymer.
$\mathscr{S}(V)$ is called the \emph{polymer of $V$}.
\end{edefin}

We must remark that there may be more than one irreducible component $W$ corresponding to the same part $T\in\mathscr{S}(V)$ of the Dynkin diagram.
We also remark that this definition using polymers, even if it is not
usual in literature where the weights of a cocharacter is used, is useful
because we want to be able to work with irreducible components of
representations over $\ql$ that are only parts of the natural representation
of the Mumford-Tate group.

\begin{etheo}
Let $A/\C$ be an abelian variety, $V=\Hun(A,\Q)$, $F=\Q$. Then the natural representation $G\to \GL(V)$ is of Hodge type. The representation $G_\ql\to\GL(V)_\ql$ obtained by scalar extension is also of Hodge type.
\end{etheo}
\begin{proof}
This theorem is a consequence of the translation by Deligne in terms of Shimura data of the work of Satake about symplectic embeddings of hermitian symetric domains. References and more details will be given in section \ref{constabdel}. For the moment, it suffices to apply the remark \ref{remhodge}.
The representation $G_\ql\to\GL(V)_\ql$ is also of Hodge type because to be 	a Hodge type representation is stable by scalars extension.
\end{proof}

\begin{erem}
The representation $G_\ql\to\GL(V)_\ql$ in the last theorem gives the action of the Mumford-Tate group on the Tate module of the abelian variety (see the proof of \ref{CRITERION}). This is the kind of representation we want to study in our combinatorial criterion.
\end{erem}

\subsection{Optimal combinatorial criterion}
As was seen, the building blocks for Hodge type representations of Mumford-Tate groups are minuscule representations and their image contains many unipotents of index $2$. We will use the galois action on tensor products to construct unipotents of higher index. The following simple lemma is central in our work.
Let $F$ be a field of characteristic $0$. 

\begin{elem}[Key lemma]
\label{keylemma}
Let $W_1$ and $W_2$ be two vector spaces over $F$ and $u_i\in~\GL(W_i)$, $i=1,2$ be two unipotents of respective index $n_1$ and $n_2$. Then $u_1\otimes u_2\in~\GL(W_1\otimes~W_2)$ is unipotent of index $n_1+n_2-1$.
\end{elem}
\begin{proof}
We use the functorial description of unipotents given in \ref{unip-reductive}. To each $u_i$ is associated a morphism $\rho_i\colon \G_a\to \GL(W_i)$ for which $\rho_i(1)=u_i$.
The image of the universal point $T=\id_{\G_a}\in \G_a(F[T])$ by $\rho_i$ is given by a polynomial $P_i$:
$$
\begin{array}{cccc}
{\rho_i}_{(F[T])}: & \G_a(F[T])	&\to    &\GL(W_i)(F[T])\\
		& T		&\mapsto&P_i=\exp(TN_i)=\sum_{j\geq 0}\frac{N_i^jT^j}{j!}
\end{array}
$$
for $N_i$ the nilpotent endomorphism given by the logarithm $$N_i=\log(\rho_i(1))=\sum_{j>0}(-1)^{j+1}(\rho_i(1)-\id)^j.$$
The degre of $P_i$ is equal to $n_i-1$.
The tensor product $u_1\otimes u_2$ corresponds to the polynomial $P_1\otimes P_2$ that is an endomorphism of
$$
(W_1\otimes_F W_2)(F[X])=(W_1\otimes_F F[X])\otimes_{F[X]} (W_2\otimes_F F[X]).
$$
The coefficient of the dominant term $X^{n_1+n_2-2}$ in $P_1\otimes P_2$ is $$\frac{N_1^{n_1-1}}{(n_1-1)!}\otimes \frac{N_2^{n_2-1}}{(n_2-1)!}$$
which is nontrivial because it is the tensor product of two nontrivial endomorphisms. This gives us a polynomial of degre $n_1+n_2-2$ and the index of unipotency of $u_1\otimes u_2$ is $n_1+n_2-2+1=n_1+n_2-1$.
\end{proof}

Let $G/F$ be a reductive group.
\begin{edefin}
\label{perftenstwisted}
Let $I$ be the set of simple factors of $\tilde{G}_{\overline{F}}$ and $\pi\colon D\to I$ the projection from the Dynkin diagram of $\tilde{G}_{\overline{F}}$ to $I$. Let $V$ be a Hodge type representation of $G$ and $\mathscr{S}(V)$ the polymer of $V$. $V$ will be called \emph{perfectly tens-twisted} if for all $i\in I$, at least one of the two following conditions is verified: 
\begin{enumerate}
\item There exists $T\in \mathscr{S}(V)$ such that $\Card(\Gal(\overline{F}/F)\cdot i\cap \pi(T))>1$,
\item The factor of $\tilde{G}$ over $F$ that contains the factor $G_i$ over $\overline{F}$ contains no nontrivial unipotent element defined over $F$.
\end{enumerate}
\end{edefin}
The first condition of the definition corresponds to the fact that the Galois group moves sufficiently the tensorial components of the representation. The second condition permits to overcome the difficulties appearing when the first condition isn't verified. We will now prove that the second condition can often be dropped. We will see in Theorem \ref{unip2hodge} a representation is perfectly tens-twisted if and only if its image does not contain any $F$-rational unipotent element of index $2$.

\begin{eprop}
We keep the notations of last definition. Suppose that $F=K_\p$ is a finite extension of $\qp$ and that no simple factor of $\tilde{G}$ is of inner type $A_n$. Then the representation $V$ is perfectly tens-twisted if and only if for all $i\in I$,
\begin{itemize}
\item there exists $T\in \mathscr{S}(V)$ such that $\Card(\Gal(\overline{K_\p}/K_\p)\cdot i\cap \pi(T))>1$.
\end{itemize}
\end{eprop}
\begin{proof}
If $I=\emptyset$, there is no condition.
Suppose that $I\neq\emptyset$. Let $H$ be any nontrivial simple factor of $\tilde{G}$ defined over $K_\p$. We will show that $H(K_\p)$ contains at least one nontrivial unipotent element. Since there is a natural isomorphism $\U(H)\cong\U(H^{ad})$ on unipotent subschemes (see \ref{prop-unip-adjoint}), it suffices to verify this for $H^{ad}$. We can write $H^{ad}=\Res_{L_\q/K_\p}H^s$ with $H^s$ absolutely simple and $L_\q/K_\p$ a finite extension. The universal property of $\Res_{L_\q/K_\p}H^s$ permits to describe its points with values in $\G_{a,K_\p}$:
$$\Res_{L_\q/K_\p}H^s(\G_{a,K_\p})\overset{\sim}{\to}H^s(\G_{a,L_\q}).$$
If we restrict to points respecting the group structure of $\G_a$, this gives a bijection
$$\U(\Res_{L_\q/K_\p}H^s)(K_\p)\overset{\sim}{\to}\U(H^s)(L_\q)$$
on unipotent subsets.
So we want to prove that $H^s(L_\q)$ contains nontrivial unipotent elements. Since $H^s$ is adjoint, this is equivalent to the fact that $H^s$ is not anisotropic by \cite{Bo5}, proof of 8.5.
But by \cite{Tits}, 3.3.3, we know that a group which is anisotropic and absolutely almost simple over a $\p$-adic field is of inner type $A_n$. By hypothesis, $H^s$ is not of this type so it is not anisotropic. This finishes the proof.
\end{proof}

\begin{etheo}
\label{unip2hodge}
Let $\rho\colon G\to \GL(V)$ be a Hodge type representation. The two following conditions are equivalent: 
\begin{enumerate}
\item $\mathrm{Im}(\rho)$ contains no nontrivial unipotents of index $2$ defined over $F$,
\item $\rho$ is perfectly tens-twisted.
\end{enumerate}
\end{etheo}

\subsection{Demonstration of theorem \ref{unip2hodge}}
Since by \ref{prop-unip-adjoint} the morphism $\tilde{G}\to G$ induces an isomorphism $\U(\tilde{G})\cong \U(G)$ on unipotent subschemes, the first condition depends only on the restriction $\tilde{G}\to \GL(V)$ of $\rho$ to $\tilde{G}$. By definition, the second condition also. So we can suppose that $G\cong \tilde{G}$. Recall that $I$ (resp. $D$) denotes the set of simple factors (resp. the Dynkin diagram) of $\tilde{G}_{\overline{F}}$ and $\mathscr{S}(V)$ the polymer of the representation.We denote $\pi\colon D\to I$ the projection.

Suppose first that $\rho$ is perfectly tens-twisted. Let $u\in G(F)$ be a nontrivial unipotent
element and let $G_i$, $i\in I$ be a quasi-simple factor of $G_{\overline{F}}$ such that the
component $u_i$ of $u$ in $G_i$ is nontrivial. We have two cases: 
\begin{itemize}
\item Suppose that the first condition for perfectly tens-twistedness is verified for $i$. There exists an irreducible component $W$ of $V_{\overline{F}}$ and a corresponding $T_W\in \mathscr{S}(V)$, such that $\Card(\Gal(\overline{F}/F)\cdot i\cap \pi(T_W))>1$. Let $j,j'$ be two points of $\Gal(\overline{F}/F)\cdot i\cap \pi(T_W)$. Then the nontrivial component $u_i$ is sent by Galois to $u_j$ and $u_{j'}$ so those two other components of $u$ are nontrivial. By definition of the polymer of a Hodge type representation, the representation $W$ is written as a tensor product of minuscule representations $\rho_k\colon G_k\to \GL(W_k)$ of the quasi-simple factors $G_k$ of $G_{\overline{F}}$ for $k\in \pi(T_W)$. Those representations are irreducible representations of quasi-simple groups, so their kernel is finite and can't contain any unipotent element. So $u_j$ and $u_{j'}$ have nontrivial image in $\GL(W_j)$ and $\GL(W_{j'})$ respectively. By the key lemma \ref{keylemma}, the index of unipotency of the tensor product of two nontrivial unipotents is strictly superior to $2$ and it can't decrease when you apply another tensor product. So the index of unipotency of $\rho(u)_{|W}=\cdots\otimes\rho_j(u_j)\otimes\cdots\otimes\rho_{j'}(u_{j'})\otimes\cdots$ in $\GL(W)$ is greater than $2$. The same is true for the image of $u$ in $\GL(V)$.
\item Suppose now that the second condition for perfectly tens-twistedness is verified for $i$. Then the factor of $G$ over $F$ that contains the component $u_i$ of $u\in G(F)$ contains no nontrivial unipotent element defined over $F$. So $u_i$ must be trivial, which is a contradiction.
\end{itemize}
We thus have proved that if $\rho$ is perfectly tens-twisted then its image contains no nontrivial unipotent of index $2$.

Suppose now that $\rho$ is not perfectly tens-twisted. By definition, there exists $i\in I$  such that
\begin{itemize}
\item \mbox{$\Card(\Gal(\overline{F}/F)\cdot i\cap \pi(T))\leq 1$} for all $T\in \mathscr{S}(V)$,
\item if we denote by $H$ the simple factor of $G$ that contains $G_i$ over $\overline{F}$, $H(F)$ contains nontrivial unipotent elements.
\end{itemize}
We recall that since the representation is of Hodge type, we can write it
$$V_{\overline{F}}=\bigoplus_{T_k\in \mathscr{S}(V)}\Big(\otimes_{s\in T_k} W_s\Big)^{\oplus n_k}$$
with $W_s$ a minuscule representation of the factor $G_{\pi(v)}$ of heighest weight $\omega_s$, $s\in T_k\subset D$ and $n_k$ the multiplicities of the irreducible components.
Let $J:=\{\Gal(\overline{F}/F)\cdot i\}\subset I$ be the set of simple factors of $H_{\overline{F}}$. By definition of a Hodge type representation, $G_i$ acts nontrivially on $V_{\overline{F}}$, so there exists $T_0\in \mathscr{S}(V)$ such that $\Card(J\cap \pi(T_0))=1$. We now write the decomposition of the part $V_{\overline{F}}'$ of $V_{\overline{F}}$ on which the action of $H$ is nontrivial.
$$V_{\overline{F}}'=
\bigoplus_{
\begin{array}{c}
\{T_k\in \mathscr{S}(V)\,|\,\Card(\pi(T_k)\cap J)=1\}\\
\{s_k\}:=T_k\cap\pi^{-1}(J)
\end{array}
}
\Bigg(W_{s_k}\otimes \bigg(
\underset{s\in T_k,s\neq s_k}{\otimes}W_s
\bigg)\Bigg)^{\oplus n_k}$$
By the existence of $T_0$, this sum is nontrivial. It gives a sub-$H$-representation $V'$ of $V$ that is defined over $F$ because it is isomorphic to $V/inv(V)$ where $inv(V)$ is the subspace of $V$ of $H$-invariant elements. The action of $H$ on $(\bigotimes_{s\in T_k,s\neq s_k}W_s)$ is trivial and the action of $H$ on each $W_{s_k}$ is nontrivial. The representation $V_{\overline{F}}'$ thus decomposes as a direct sum of minuscule representations of the simple factors of $H_{\overline{F}}$. Let $T\subset H$ be a maximal torus and $R=R(H,T)$ be the corresponding scheme of roots. Since $H(F)$ contains nontrivial unipotent elements, the structure theory of reductive groups imply that $R(F)$ contains at least one element $r$. Indeed, by \cite{Bo5}, 8.3, we know that any nontrivial unipotent $F$-subgroup is contained in a nontrivial $F$-parabolic subgroup  and by \cite{Bo6}, 6.5, (5), (6), we know that any nontrivial $F$-parabolic subgroup is conjugated over $F$ to a standard parabolic subgroup that is naturally associated to a nontrivial subset $\Theta\subset R(F)$, so $R(F)$ can't be trivial. Over $\overline{F}$, the root $r$ decomposes as a sum $r=\sum_{j\in J}r_j$ of roots of the simple factors $(G_j,T\cap G_j), j\in J$ of $(H_{\overline{F}},T_{\overline{F}})$.  The corresponding morphism $u_r\colon \G_a\to H$ also decomposes as a sum $u_r=\sum_{j\in J}u_{r_j}=\G_a\to \prod_{j\in J}G_j$ of nontrivial morphisms. So the image of this morphism in $V_{\overline{F}}'$ is a direct sum of images of standard unipotent subgroups of the $G_j,j\in J$ in minuscule representations of those factors. It thus must be unipotent of index $2$, by definition \ref{definminuscule} of minuscule representations. \hfill$\Box$

\subsection{Application to the conjecture}
We now have a combinatorial criterion for potentially good reduction.
Let $K$ be a number field and $A/K$ an abelian variety. Let $G=\MT(A)$ and $\rho\colon G\to \GL(V)$ be the natural representation on $V=\Hun(A,\Q)$. 
\begin{etheo}
\label{combincriteria}
If there exists a prime $\ell$ of $\Q$ such that $\rho_\ql\colon G_\ql\to \GL(V)_\ql$ is perfectly tens-twisted, then $A$ has potentially good reduction at every finite place of $K$.
\end{etheo}
\begin{proof}
This is a direct consequence of our theorical criterion \ref{CRITERION} for potentially good reduction and of theorem \ref{unip2hodge}.
\end{proof}

Suppose for the rest of this paragraph that $A$ is a simple abelian variety with
Mumford-Tate group $G$ such that $G^{ad}$ is non nontrivial $\Q$-simple.
Let $I$ (resp. $D$) be the nontrivial set of simple factors (resp. the Dynkin diagram) of $\tilde{G}_{\qb}$ and $\mathscr{S}(V)$ be the polymer of the representation $V$ of $G$. We denote $\pi\colon D\to I$ the natural projection. We now give a simple condition for the existence of a prime $\ell$ such that $V_\ql$ is perfectly tens-twisted. This is the criterion we usually apply to obtain our new results.
\begin{eprop}
\label{combinappli}
Suppose that there exists $\sigma\in\Gal(\qb/\Q)$ such that for all $i\in I$, there exists $T\in \mathscr{S}(V)$ such that
$$\Card(\langle\sigma\rangle\cdot i\cap \pi(T))>1$$
then there exists a prime $\ell$ such that $V_\ql$ is perfectly tens-twisted and $A$ has potentially good reduction at every discrete place of $K$.
\end{eprop}
\begin{proof}
By Tchebotarev's theorem, there exist a prime $\ell$ of $\Q$ and a prime $\bar{\ell}$ of $\qb$ over $\ell$, such that
$$\sigma\in \DD(\bar{\ell}),$$
where $\DD(\bar{\ell})\subset \Gal(\qlb/\ql)$ denotes the decomposition group at $\bar{\ell}$.
The prime $\qb\overset{\bar{\ell}}{\to} \qlb$ gives an isomorphism $\Gal(\qlb/\ql)\cong \DD(\bar{\ell})$.
We denote by $I$ (resp. $I_\ell$) the set of simple factors of $\tilde{G}_\qb$ (resp. $\tilde{G}_\qlb$). The embedding $\qb\overset{\bar{\ell}}{\to}\qlb$ gives a $\Gal(\qlb/\ql)$-equivariant bijection $I\overset{\bar{\ell}}{\to}I_\ell$. If $i\in I_\ell$, we have
$$\Card(\Gal(\qlb/\ql)\cdot i\cap \pi(T))\geq \Card(\langle\sigma\rangle\cdot i\cap \pi(T))>1$$
so the representation $V_\ql$ is perfectly tens-twisted. Theorem \ref{combincriteria} gives the final result.
\end{proof}

\begin{ecor}[Cyclic case]
\label{corcyclic}
If the action of $\Gal(\qb/\Q)$ on $I$ is cyclic and if there exists $T\in \mathscr{S}(V)$ such that $\Card(\pi(T))>1$ then there exists a prime $\ell$ of $\Q$ such that $V_\ql$ is perfectly tens-twisted. In that case, $A$ has potentially good reduction at every discrete place of $K$.
\end{ecor}
\begin{proof}
By hypothesis, the action of Galois on $I$ is transitive.
The generator $\sigma$ of the action of Galois on $I$ verifies that for all $i\in I$, $\Card(\langle\sigma\rangle\cdot i\cap \pi(T))>1$. The proposition \ref{combinappli} implies that $A$ has potentially good reduction at every discrete place of $K$.
\end{proof}

\begin{ecor}[General Mumford type]
\label{cormumford}
If $\Card(I)>1$ and if there exists $T\in~\mathscr{S}(V)$ such that $\pi(T)=I$ then there exists a prime $\ell$ of $\Q$ such that $V_\ql$ is perfectly tens-twisted and $A$ has potentially good reduction at every discrete place of $K$.
\end{ecor}
\begin{proof}
The action of the Galois group on $I$ is transitive.
The following simple lemma on finite group actions gives a $\sigma\in \Gal(\qb/\Q)$ that acts on $I$ without fixed points. This $\sigma$ verifies the condition that  for all $i\in I$, $\Card(\langle\sigma\rangle\cdot i\cap \pi(T))>1$, which is the hypothesis of \ref{combinappli}.
\end{proof}

\begin{elem}
Let $n>1$ be an integer and $G\subset S_n$ a subgroup that acts transitively on $E=\{1,...,n\}$. Then there exists an element $g\in G$ without fixed point.
\end{elem}
\begin{proof}
If $g$ is an element of $G$, we denote $F(g)$ the number of fixed points of $g$ in $E=\{1,...,n\}$. We then have the formula
$$\sum_{g\in G}F(g)=\Card(G).$$
Indeed, denote $S=\{(a,g)\in E\times G\,|\,ga=a\}$.
We will count points of $S$ in two different maners. For all $a\in E$, the number of $g$ with $(a,g)\in S$ is the order of the stabilizer $G_a$ of $a$. Since the action is transitive, all the stabilizers are of the same order and we have $\Card(S)=|G_a|.|E|=|G|$. Now, for all $g\in G$, the number of $(a,g)\in S$ is $F(g)$. So we obtain $\sum_{g\in G}F(g)=|S|=|G|$.
Now suppose that every nontrivial element of $G$ has at least one fixed point. This gives us the inequality
$$|G|-n=\sum_{g\in G-{\id}}F(g)\geq |G|-1$$
which is impossible since $n>1$.
\end{proof}

\section{Existence of abelian varieties satisfying our criterion}
\label{constab}
We will now construct, using the machinery of Shimura varieties, abelian varieties satisfying our criterion. The Shimura varieties will be constructed from adjoint Shimura data $(G,X)$ and some particular polymers for $G$ depending on  $X$, called Deligne polymers. Our abelian varieties will be fibres of universal families of abelian varieties over those Shimura varieties.

Since our combinatorial criterion focuses on the properties of the polymer associated to the natural representation of the Mumford-Tate group, this gives us a concrete way to construct a lot of examples of abelian varieties satisfying our criterion.

\subsection{Shimura varieties and Morita's conjecture}
A \emph{Shimura datum} is a pair $(G,X)$, with $G$ a reductive algebraic group over $ \Q$ and $X$ a $G(\R)$-conjugacy class in the set of morphisms of algebraic groups $\Hom(\mathbb{S},G_{\R})$ satisfying the three conditions of Deligne \cite{De4}, 2.1.1.1-3. These conditions imply that $X$ is a hermitian symetric domain, such that every representation of $G$ on a $\Q$-vector space defines a polarizable variation of Hodge structure on $X$. For $(G,X)$ a Shimura datum and $K$ a compact open subgroup of $G(\A_f)$, we let $\Sh_K(G,X)(\C)$ denote the complex analytic variety $G(\Q)\backslash (X\times G(\A_f)/K)$, which has a natural structure of quasi-projective complex variety, denoted $\Sh_K(G,X)_\C$; we denote $\Sh(G,X)_\C$ the projective limit over all $K$. A morphism of Shimura data from $(G_1,X_1)$ to $(G_2,X_2)$ is a morphism $f\colon G_1\to G_2$ that maps $X_1$ to $X_2$; such an $f$ induces a morphism $\Sh(f)$ from $\Sh(G_1,X_1)_\C$ to $\Sh(G_2,X_2)_\C$. To every point $x$ in $\Sh(G_1,X_1)(\C)$ is associated a \emph{Mumford-Tate group} $\MT(x)$ that is the smallest $\Q$-algebraic subgroup of $G_1$ containing the image of $h_x\colon \mathbb{S}\to G_{1,\R}$ over $\R$. The Shimura variety $\Sh(G_1,X_1)_\C$ contains some particular points of great interest for the study of its arithmetic properties, called \emph{special points}. These are the points $x$ of $\Sh(G_1,X_1)_\C$ the Mumford-Tate group of which is a torus.

Now take an abelian variety $A$ over a number field $F$. Suppose given an embedding $F\hookrightarrow \C$ and let $\MT(A)$ denote the Mumford-Tate group of $A_\C$. We can associate to $A$ a Shimura variety that gives roughtly speaking the ``biggest family of abelian varieties with Mumford-Tate groups contained in $\MT(A)$'' in the following way. Let $h\colon \mathbb{S}\to \MT(A)_\R$ the natural morphism. Let $X_A$ be the $\MT(A)(\R)$-conjugacy class of $h$ in $\Hom(\mathbb{S},\MT(A)_\R)$. Then $(\MT(A),X_A)$ is a Shimura datum. Let $V=\Hun(A_\C,\Q)$ be the Hodge structure of $A$. It is naturally a polarizable Hodge structure. The choice of a polarization gives a symplectic form $\psi$ on $V$ such that the natural morphism $\MT(A)\to \GL(V)$ factors through $\GSp(V,\psi)$, which we denote more simply $\GSp(V)$, with $\psi$ understood. We denote by $S^\pm$ the double Siegel Half space corresponding to $(V_\R,\psi_\R)$ (see Deligne \cite{De4}, 1.3.1). $(\GSp(V),S^\pm)$ is also a Shimura datum and we have an embedding of Shimura data
$$(\MT(A),X_A)\hookrightarrow (\GSp(V),S^\pm)$$
that is called a \emph{Siegel embedding of the Shimura datum $(\MT(A),X_A)$}.
More generally, an embedding
$$(G_1,X_1)\hookrightarrow (\GSp(V),S^\pm)$$
of a Shimura datum $(G_1,X_1)$ is called a \emph{Siegel embedding} and the corresponding Shimura variety is called of \emph{Hodge type}.

Shimura varieties admit models over number fields, called canonical models, enjoying nice properties. Take a Shimura datum $(G_1,X_1)$ and denote $E=E(G_1,X_1)$ the reflex field. Then $\Sh(G_1,X_1)_\C$ admits a canonical model $\Sh(G_1,X_1)_E$ in the sense of \cite{De4}, 2.7.10. Roughly speaking, this is a model such that special points are algebraic and the Galois action on them is given in an explicit way. The canonical model of the Siegel modular variety $\Sh(\GSp(V),S^\pm)$ is defined over $\Q$ and if we have a Siegel embedding of a Hodge type Shimura datum
$(G_1,X_1)\hookrightarrow~(\GSp(V),S^\pm)$,
it induces by functoriality of canonical models a morphism
$\Sh(G_1,X_1)_E\hookrightarrow~\Sh(\GSp(V),S^\pm)_E$. After the choice of a level structure, the right term has an interpretation as a fine moduli space for principally polarized abelian varieties over $E$. So if we choose sufficiently small compact open subgroups $K_1\subset G_1(\A_f)$ and $K_S\subset \GSp(V)(\A_f)$, we can suppose that $\Sh_{K_S}(\GSp(V),S^\pm)_E$ carries a universal family $\mathcal{A}_u$ of abelian varieties that pulls back over $\Sh=\Sh_{K_1}(G_1,X_1)_E$, giving an abelian scheme $\mathcal{A}$ that fits into the following cartesian diagram:
$$
\xymatrix{
\mathcal{A}\ar[r]\ar[d] & \mathcal{A}_u\ar[d]\\
\Sh      \ar[r]       &  \Sh_{K_S}(\GSp(V),S^\pm)_E
}
$$

\begin{edefin}
\label{absch-siegel}
Let $(G_1,X_1)\hookrightarrow~(\GSp(V),S^\pm)$ be a Siegel embedding. After the choice of a sufficiently small compact open subgroup $K_1\subset G_1(\A_f)$, the above construction gives an abelian scheme $\mathcal{A}\to \Sh$ over the canonical model $\Sh=\Sh_{K_1}(G_1,X_1)_E$ of the corresponding Shimura variety.
This abelian scheme is called the \emph{family of abelian varieties over} $\Sh$ associated to the Siegel embedding $(G_1,X_1)\hookrightarrow~(\GSp(V),S^\pm)$.
\end{edefin}

\begin{erem}
\label{remabelianfibre}
There is an interpretation of canonical models of Hodge type Shimura varieties as moduli spaces for principally polarized abelian varieties with absolute Hodge cycles (see \cite{Brylinski}, 1.2.6 for an account of this idea).
It implies that if we take an abelian variety $A$ over a number field $F$, and $(\MT(A),X_A)$ the associated Shimura datum, there exists a finite extension of the composite field $F.E(\MT(A),X_A)$ such that $A$ is a fibre of the natural family over $\Sh_K(\MT(A),X_A)$ for an appropriate compact open subgroup $K\subset \MT(A)(\A_f)$. This shows that \emph{every abelian variety over a number field lives in the family of abelian varieties associated to a Siegel embedding of a Hodge type Shimura datum}.
\end{erem}

The original formulation of Morita's conjecture was made by Morita in terms of {\PEL} Shimura varieties. We reformulate it for Hodge type Shimura varieties.
\begin{econj}[Morita]
\label{conjgeomorita}
Suppose given a Siegel embedding of a Hodge type Shimura datum $$(G_1,X_1)\hookrightarrow~(\GSp(V),S^\pm)$$
and a compact open subgroup $K_1\subset G_1(\A_f)$.
If $\Sh_{K_1}(G_1,X_1)_\C^{an}$ is a compact topological space then the fibre of the family $\mathcal{A}\to~\Sh$ associated to the prescribed Siegel embedding (see \ref{absch-siegel}) at any point of $\Sh=\Sh_{K_1}(G_1,X_1)_E$ over a number field $F$ is an abelian variety $A$ with potentially good reduction at every discrete place of $F$.
\end{econj}

\begin{eprop}
Our formulation \ref{conjalgmorita} of Morita's conjecture is equivalent to the original formulation \ref{conjgeomorita}.
\end{eprop}
\begin{proof}
Suppose first that \ref{conjgeomorita} is true. Let $A$ be an abelian variety over a number field $F$ which Mumford-Tate group $\MT(A)$ contains no nontrivial unipotent element over $\Q$. Let $(G_1,X_1)$ be the Shimura datum $(\MT(A),X_A)$ associated to $A$. The lemma \ref{lem-compacite} shows that $\Sh=\Sh_{K_1}(G_1,X_1)_\C$ is compact. The universal family defined in \ref{absch-siegel} gives an abelian scheme $\mathcal{A}\to \Sh$ corresponding to the natural Siegel embedding of $(G_1,X_1)=(\MT(A),X_A)$.
Remark \ref{remabelianfibre} shows that there exists a point $x\in \Sh(F')$ in a number field $F'\supset F$ such that $\mathcal{A}_x=A_{F'}$. Our formulation \ref{conjalgmorita} of Morita's conjecture says  that $A_{F'}$ has potentially good reduction at every discrete place of $F'$. This proves one implication.

Suppose now that \ref{conjalgmorita} is true. Let $(G_1,X_1)$ be a Hodge type Shimura datum, $K_1\subset G_1(\A_f)$ be a compact open subgroup, and suppose that $\Sh_{K_1}(G_1,X_1)_\C^{an}$ is compact. Lemma \ref{lem-compacite} tells us that $G_1(\Q)$ contains no nontrivial unipotent element. The Mumford-Tate group of any fibre of the universal family $\mathcal{A}\to~\Sh$ of \ref{absch-siegel} is contained in $G_1$ and thus contains no nontrivial unipotent element. Any fibre of this abelian scheme verifies the hypothesis of our formulation \ref{conjalgmorita} of Morita's conjecture and has potentially good reduction at every discrete place of its field of definition. This proves the equivalence of the two formulations.
\end{proof}

\begin{elem}
\label{lem-compacite}
Let $(G,X)$ be a Shimura datum and $K\subset G(\A_f)$ be a compact open subgroup.
The Shimura variety $\Sh_K(G,X)(\C)$ is compact if and only if $G(\Q)$ contains no nontrivial unipotent element.
\end{elem}
\begin{proof}
Along this proof, the exponent $^+$ will denote a topological connected component.
Let $\pi\colon G\to G^{ad}$ be the natural projection and $K^{ad}\subset G^{ad}(\A_f)$ be a compact open subgroup that contains $\pi(K)$. Let $X^+$ be a connected component of $X$. The projection $X\to X^{ad}$ induces an isomorphism from $X^+$ to a connected component $X^{ad+}$ of $X^{ad}$.
We still denote $\pi$ the natural morphism
$\Sh_K(G,X)(\C)\to \Sh_{K^{ad}}(G^{ad},X^{ad})(\C).$
Let $\Sh_{K^{ad}}(G^{ad},X^{ad})(\C)^+$ be a connected component of the adjoint Shimura variety and $\Sh_K(G,X)(\C)^+$ be a connected component of the Shimura variety the image of which is contained in $\Sh_{K^{ad}}(G^{ad},X^{ad})(\C)^+$.
By \cite{De4}, 2.1.2, we know that the morphism
$$\Sh_K(G,X)(\C)^+\to \Sh_{K^{ad}}(G^{ad},X^{ad})(\C)^+$$
is given by a projection
$$\Gamma\backslash X^+\to \Gamma'\backslash X^+$$
with $\Gamma$ and $\Gamma'$ two arithmetic subgroups of the connected component of identity $G^{ad}(\R)^+$ of $G^{ad}(\R)$. This quotient is a finite morphism and we know (see \cite{De4}, 2.1.2) that the Shimura varieties in play have finitely many connected components so the natural morphism
$\pi\colon \Sh_K(G,X)(\C)\to \Sh_{K^{ad}}(G^{ad},X^{ad})(\C)$
is finite.
Thus $\pi$ is proper, so the Shimura variety is compact if and only if its adjoint Shimura variety is compact.
Since the center of $G^{ad}$ is trivial, the conditions of \cite{De4}, 2.1.10 are verified and we know that
$$\Sh(G^{ad},X^{ad})(\C)=G^{ad}(\Q)\backslash X\times G^{ad}(\A_f).$$
We also know from \cite{De4}, 2.1.4 that $G(\A_f)$ acts on $\Sh(G^{ad},X^{ad})$ and that $\Sh(G^{ad},X^{ad})/K^{ad}=\Sh_{K^{ad}}(G^{ad},X^{ad})$.
Since $K^{ad}$ is compact and acts continuously on the separated topological space $\Sh(G^{ad},X^{ad})(\C)$, the projection map
$$\Sh(G^{ad},X^{ad})\to \Sh_{K^{ad}}(G^{ad},X^{ad})$$
is proper by \cite{Boutopo}, Chap 3, \S4, n$^\circ$ 1, Prop 2 so $\Sh_{K^{ad}}(G^{ad},X^{ad})(\C)$ is compact if and only if $\Sh(G^{ad},X^{ad})(\C)$ is compact.
Writing $X^{ad}=G^{ad}(\R)/K^{ad}_\infty$ with $K^{ad}_\infty$ a maximal compact subgroup of $G^{ad}(\R)$, we obtain a continuous quotient map
$$G^{ad}(\Q)\backslash G^{ad}(\A)\to G^{ad}(\Q)\backslash G^{ad}(\R)\times G^{ad}(\A_f)/K^{ad}_\infty\cong \Sh(G^{ad},X^{ad})(\C)$$
which is proper because $G^{ad}(\Q)\backslash G^{ad}(\A)$ is separated and $K^{ad}_\infty$ is compact \cite{Boutopo}, Chap 3, \S4, n$^\circ$ 1, Prop 2.
By a result of Borel and Harish-Chandra (see \cite{Bo4}, 5.6, p21), the compacity of $G^{ad}(\Q)\backslash G^{ad}(\A)$ is equivalent to the fact that $G^{ad}(\Q)$ contains no nontrivial unipotent element. So the Shimura variety $\Sh(G,X)(\C)$ is compact if and only if $G^{ad}(\Q)$ contains no nontrivial unipotent element. The isomorphism $\U(G)\cong\U(G^{ad})$ of unipotent subschemes given in \ref{prop-unip-adjoint} shows that this is also equivalent to the fact that $G(\Q)$ contains no nontrivial unipotent element.
\end{proof}

\subsection{Construction of families of abelian varieties using Deligne polymers}
\label{constabdel}
let $(G,X)$ be a Shimura datum with $G$ adjoint simple. Then Deligne \cite{De4}, 1.3,2.3.7 using Satake's work \cite{Satake3}, classifies diagrams
$$(G,X)\leftarrow (G_1,X_1) \hookrightarrow (\GSp(V),S^\pm)$$
with $(G_1,X_1)^{ad}=(G,X)$.

By \cite{De4}, 2.3.4(a) we know that $G=\Res_{F/\Q}G^s$ with $F$ totally real and $G^s$ absolutely simple over $F$.
Let $I$ be the set of factors of $G_{\R}$ (or, which is the same, the set of real embeddings of $F$) and $I_c$ (resp. $I_{nc}$) be the subset of indices corresponding to compact (resp. non compact) factors of $G_{\R}$. Let $D$ be the Dynkin diagram of $G$, $D_{c}$ (resp. $D_{nc}$) the part corresponding to compact (resp. non compact) factors of $G_{\R}$. If $i\in I$, we denote $D_i$ the Dynkin diagram of the corresponding factor. The action of $\Gal(\qb/\Q)$ is compatible with the projection $\pi\colon D\to I$. Let $\mathscr{P}(D)$ be the set of parts of the set of vertices of the Dynkin diagram. $X$ gives a special vertex $s_i$ of $D_i$ for every $i\in I_{nc}$ (This vertex is given in \cite{De4}, 1.2.5 and the set of all those special vertices for $i\in I_{nc}$ uniquely determines $X$).

\begin{edefin}[Deligne polymer]
\label{defindelignepolymer}
A polymer $\mathscr{S}\subset \mathscr{P}(D)$ is called a \emph{Deligne polymer} for $(G,X)$ if
\begin{itemize}
\item for all $T\in \mathscr{S}$, $T\cap D_{nc}$ is empty or reduced to one point $s_T\in D_i$ ($i\in I_{nc}$), and in the table \cite{De4}, 1.3.9 for $(D_i,s_i)$, $s_T$ is an underlined vertex.
\end{itemize}
\end{edefin}
Let
$$(G,X)\leftarrow (G_1,X_1) \hookrightarrow (\GSp(V),S^\pm)$$
be a diagram with $(G_1,X_1)^{ad}=(G,X)$ and let $\tilde{G}$ be the universal covering of $G$. The natural map $\tilde{G}\to G$ lifts to $\tilde{G}\to G_1$ and we obtain a representation $V_\C$ of $\tilde{G}$. The Satake classification tells us that every irreducible component $W$ of $V_\C$ is of the form $\otimes_{i\in I_W} W_i$ for some $I_W\subset I$, with $W_i$ a fundamental representation of $G_i$ corresponding to a vextex $\tau(i)$ of $D_i$. Denote $T_W$ the subset $\tau(I_W)\subset D$ of the Dynkin diagram. We thus know that $W$ is of the form $\otimes_{s\in T_W} W_s$ for $W_s$ a fundamental representation corresponding to the vertex $s\in T_W\subset D$ of the Dynkin diagram. If we denote $\mathscr{S}(V)$ the set of $T_W\subset D$ for $W\subset V_\C$ irreducible, then $\mathscr{S}(V)$ is a Deligne polymer.

\begin{edefin}
Let $$(G,X)\leftarrow (G_1,X_1) \hookrightarrow (\GSp(V),S^\pm)$$ be a diagram with $(G_1,X_1)^{ad}=(G,X)$. The Deligne polymer obtained in the preceding description is called the \emph{Deligne polymer for $(G,X)$ associated to the Siegel embedding $(G_1,X_1) \hookrightarrow (\GSp(V),S^\pm)$}.
\end{edefin}

\begin{erem}
\label{remhodge}
The representation $G_1\to \GL(V)$ is of Hodge type over $\Q$ (see Definition \ref{rephodgetype}) and the Deligne polymer $\mathscr{S}(V)$ is a polymer for $G$ in the sense given in \ref{defin-polymer} with an additional condition on vertices that depends on $X$. In particular, let $A$ be an abelian variety over a number field $F\hookrightarrow\C$, let $V$ be $\Hun(A_\C,\Q)$ and
$$(\MT(A),X_A)\hookrightarrow (\GSp(V),S^\pm)$$
be the corresponding Siegel embedding.
Suppose that $\MT(A)$ has no commutative simple factor. We obtain that the natural representation
$$\MT(A)\to \GL(V)$$
is of Hodge type with polymer $\mathscr{S}(V)$ equal to the Deligne polymer associated to this Siegel embedding.
\end{erem}

\begin{eprop}
\label{constab-polymer}
Let $(G,X)$ be a Shimura datum with $G$ a $\Q$-simple group of adjoint type. For every Deligne polymer $\mathscr{S}$ for $(G,X)$, there exists a (not necessarily unique) Siegel embedding of a Hodge type Shimura datum $(G_1,X_1)\hookrightarrow (\GSp(V),S^\pm)$ with adjoint Shimura datum $(G,X)$ and with a polymer equal to $\mathscr{S}$. After the choice of a convenient compact open subgroup $K_1\subset G_1(\A_f)$, this gives an algebraic variety $\Sh=\Sh_{K_1}(G_1,X_1)_E$ over a number field $E$ and an abelian scheme $\mathcal{A}\to \Sh$ for which every fibre $\mathcal{A}_x$ over a point $x\in \Sh(E')$ in a finite extension of $E$ gives an abelian variety with Mumford-Tate representation $$\MT(\mathcal{A}_x)\subset G_1\to \GL(V)$$ contained in the Hodge type representation $G_1\to \GL(V)$ with polymer $\mathscr{S}$.
\end{eprop}
\begin{proof}
The construction of a Siegel embedding 
$(G_1,X_1)\hookrightarrow (\GSp(V),S^\pm)$
associated to a Deligne polymer is given in \cite{De4}, 2.3.10.
The family associated to the Siegel embedding (see \ref{absch-siegel}) provides us, after the choice of a sufficiently small compact open subgroup $K_1\subset G_1(\A_f)$, an abelian scheme $p\colon \mathcal{A}\to \Sh$ over the canonical model $\Sh=\Sh_{K_1}(G_1,X_1)_E$ with $E=E(G_1,X_1)$ that is the pull-back of the universal abelian scheme $\mathcal{A}_u$ on the Siegel modular variety $\Sh_K(\GSp(V),S^\pm)_E$ for a sufficiently small compact open subgroup $K\subset \GSp(V)(\A_f)$. So we have the following cartesian diagram:
$$
\xymatrix{
\mathcal{A}_\C\ar[r]\ar[d]^p & \mathcal{A}_{u,\C}\ar[d]_{p_u}\\
\Sh_\C      \ar[r]     	  &  \Sh_{K_S}(\GSp(V),S^\pm)_\C
}
$$
By definition, the variation of $\Q$-Hodge structure $\mathbb{R}^1p_{u,*}\Q$ on the Siegel modular variety is given by the representation $\GSp(V)\to\GL(V)$  (in the sense of \cite{Moonen1}, 2.3), i.e. by the trivial bundle $S^\pm\times V$ with variable filtration on $V_\C$ depending on points in $S^\pm$ (see also \cite{De6}, 4).
Since $\mathcal{A}_\C$ is the pull-back $\mathcal{A}_{u,\C}$, the variation of $\Q$-Hodge structures $\mathbb{R}^1p_*\Q$ is the pull back of $\mathbb{R}^1p_{u,*}\Q$ and it is given (in the sense of \cite{Moonen1}, 2.3) by the representation $G_1\to\GL(V)$. 
So the representation of the Mumford-Tate group $\MT(\mathcal{A}_x)\to \GL(\mathbb{R}^1p_*\Q)\cong\GL(V)$ for every point $x$ of $\Sh_\C$ factors through the Hodge type representation $G_1\to \GL(V)$ (see \cite{Moonen1}, 2.3 for the factorisation and \cite{Moonen1}, 1.2 for facts about generic Mumford-Tate groups).
\end{proof}

\begin{edefin}[Abelian varieties associated to a Deligne polymer]
For notation simplicity, the fibres over points in number fields of the family $\mathcal{A}\to \Sh$ constructed in \ref{constab-polymer} will be called \emph{abelian varieties associated to the triplet $(G,X,\mathscr{S})$}.
\end{edefin}

\begin{eexm}
\label{exmabfibre}
If $A$ is a simple abelian variety defined over a number field $F$ with noncommutative Mumford-Tate group, $(G,X)=(\MT(A),X_A)^{ad}$, $V=\Hun(A_\C,\Q)$ and $\mathscr{S}$ is the Deligne polymer of the Siegel embedding
$$(\MT(A),X_A)\hookrightarrow (\GSp(V),S^\pm),$$
then \ref{remabelianfibre} shows that $A$ is an abelian variety associated to the triplet $(G,X,\mathscr{S})$.
\end{eexm}

\begin{eprop}
\label{constructive-result}
Let $(G,X)$ be a classical type Shimura datum with $G$ a $\Q$-simple group of adjoint type. Suppose that there exists a Deligne polymer $\mathscr{S}$ for $(G,X)$ such that the corresponding representation is perfectly tens-twisted over $\ql$ for some prime $\ell$ of $\Q$. Then the fibers of an abelian scheme $\mathcal{A}\to \Sh$ corresponding to this polymer by \ref{constab-polymer} at points of $\Sh$ in number fields have potentially good reduction at every discrete place of their field of definition.
\end{eprop}
\begin{proof}
We keep the notations of \ref{constab-polymer}. Let $x\in \Sh(E')$ be a point of $\Sh$ in a number field $E'$. The corresponding fibre $\mathcal{A}_x$ has a Mumford-Tate representation
$\MT(\mathcal{A}_x)\subset G_1\to \GL(V)$
that factors through the Hodge type representation $G_1\to \GL(V)$ with polymer $\mathscr{S}$. Since this representation is perfectly tens-twisted over $\ql$, its image contains no nontrivial unipotent of index $2$ over $\ql$, and the same is true for the image of the representation $\MT(\mathcal{A}_x)\to \GL(V)$. Using our criterion \ref{CRITERION}, we obtain that $\mathcal{A}_x$ has potentially good reduction at every discrete place of $E'$.
\end{proof}

Examples of Deligne polymers fulfiling the hypothesis of this proposition will be given in
\ref{explicit-examples}. The question of existence of such Deligne polymers is essentially
a combinatorial question about finite group actions. We must remark that it is
easy to construct a compact Shimura variety for which there exist no such
Deligne polymer, for example using the multiplicative group of a quaternion algebra $D$
(over a totally real number field) which splits at every real place but does not split at some
finite place.

\subsection{Our cases are new}
\label{PEL}
In this section, we show that the cases we are studying are almost all non {\PEL} cases of Shimura data, so are new cases of the conjecture.

Let $L$ be a semi-simple algebra of finite dimension over $\Q$ with an involution $\sigma_L$, and $V$ a finite dimensional vector space over $\Q$ with a faithful $L$-module structure and an antisymmetric non degenerated form $\psi$ such that $\psi(ax,y)=\psi(x,\sigma_L(a)y)$ for all $a\in L$ and all $x,y\in V$. This form permits to extend the involution $\sigma_L$ to $\End(V)$: if $f\in\End(V)$, we define $\sigma(f)$ by $\psi(v,\sigma(f)(w))=\psi(f(v),w)$ for all $v,w\in V$. 
We furthermore suppose that $L$ and $V$ are free over the center $Z$ of $L$.
Let
$\GSp_{L,\G_{m,\Q}}(V,\psi)$ be the group of $\G_{m,\Q}$-symplectic similitudes of $V$ given for every commutative $\Q$-algebra $A$ by
$$
\GSp_{L,\G_{m,\Q}}(V,\psi)(A):=\{a\in \GL_L(V)(A)\,|\,\sigma(a)a\in \G_{m,\Q}(A)\}.
$$
Let $G$ be the identity component of $\GSp_{L,\G_{m,\Q}}(V,\psi)$.
Assume there is a morphism $h\colon \mathbb{S}\to G_\R$ such that the Hodge structure $(V,h)$ is of type $\{(-1,0),(0,-1)\}$ and $2\pi i\psi$ is a polarization for $(V,h)$. Then the involution $\sigma$ is positive and when we take $X$ to be the set of $G(\R)$ conjugates of $h$, $(G,X)$ satisfies the axioms of a Shimura variety and it has a natural embedding in the Siegel Shimura datum $(\GSp(V,\psi),S^\pm)$. The group $G$ is determined by $L$ and $\psi$.

\begin{edefin}[{\PEL} type Shimura datum]
A Shimura datum arising from the above construction is called of \emph{{\PEL} type}.
\end{edefin}

The following proposition follows directly from the work of Kottwitz in \cite{Kottwitz}.
\begin{eprop}
\label{pel-minuscule}
Let $(G,X)\hookrightarrow (\GSp(V),S^\pm)$ be the Siegel embedding of a {\PEL} type Shimura datum and suppose that $G$ is quasi-simple and non-commutative. Then the corresponding representation $G\to \GL(V)$ decomposes over $\qb$ in a direct sum of minuscule representations of the quasi-simple factors of $\tilde{G}_{\qb}$. The corresponding polymer $\mathscr{S}(V)$ is composed of one element parts of the Dynkin diagram and this Hodge type representation contains no nontrivial tensor products of representations of the simple factors of $\tilde{G}_\qb$.
\end{eprop}

\begin{ecor}
\label{Bn-paspel}
Let $(G,X)$ be a Shimura datum with $G$ of $B_n$ type. Then $(G,X)$ is not of {\PEL} type and there exists no Shimura datum $(G_1,X_1)$ of {\PEL} type such that $(G_1,X_1)^{ad}=(G,X)$.
\end{ecor}
\begin{proof}
Since any Shimura datum $(G_1,X_1)$ with adjoint Shimura datum $(G,X)$ is also of $B_n$ type, it suffices to prove that if $(G,X)$ is of $B_n$ type then it is not of {\PEL} type. This is well known (see \cite{Kottwitz}).
\end{proof}

\section{Morita's conjecture is independent of the center}
\label{indepcenter}
\subsection{General results}
Let $F$ be a number field with an embedding $F\subset \C$ and $A$ be an abelian variety over $F$.
Recall that Morita's conjecture \ref{conjalgmorita} for $A$ states that if $\MT(A)(\Q)$ contains no nontrivial unipotent element then $A$ has potentially good reduction at every discrete place of $F$. The hypothesis of this conjecture depends only on $\MT(A)^{ad}$ because of the isomorphism $\U(\MT(A))\cong \U(\MT(A)^{ad})$ between unipotents subschemes given in \ref{prop-unip-adjoint}. We will show that the results we obtain about the conjecture are also almost only dependent of the adjoint group $\MT(A)^{ad}$.
A more precise statement is now given.

We use the category of absolute Hodge motives over $F$ constructed in \cite{De1}.
\begin{edefin}
Let $\p$ be a finite prime of $F$, $\bar{\p}$ be a prime of $\overline{F}$ over $\p$.
An absolute Hodge motive $M$ over $F$ will be called \emph{unramified at $\p$} if there exists some prime $\ell$ of $\Q$ such that the $\ell$-adic realization $M_\ell$ of $M$ is an unramified galois representation at $\p$, i.e. the inertia group $\II(\bar{p})$ acts trivially on $M_\ell$. $M$ will be called \emph{potentially unramified at $\p$} if there exists a finite extension $F'$ of $F$ and a finite prime $\p'$ of $F'$ over $\p$ such that $M_{F'}$ is unramified at $\p'$. 
\end{edefin}

\begin{eprop}
\label{adjointtip}
Let $F$ be a number field endowed with an embedding $F\subset \C$, and $A_1$, $A_2$ be two simple abelian varieties over $F$ with non-commutative Mumford-Tate group, and suppose that
\begin{itemize}
\item $\MT(A_1)^{ad}=\MT(A_2)^{ad}$,
\item the projections $h_i\colon \mathbb{S}\to \MT(A_i)_\R\to\MT(A_i)^{ad}_\R$, $i=1,2$ of the morphisms given by the Hodge structures are equal.
\end{itemize}
Then there exists an abelian absolute Hodge motive $M$ over a finite extension of $F$ with Mumford-Tate group $\MT(A_i)^{ad}$, that is a subquotient of a tensor construction on the absolute Hodge $\hunah$ of $A_1$ and $A_2$, and such that, for any finite prime $\p$ of $F$, the following statements are equivalent:
\begin{itemize}
\item $M$ is potentially unramified at $\p$,
\item $A_1$ has potentially good reduction at $\p$,
\item $A_2$ has potentially good reduction at $\p$.
\end{itemize}
\end{eprop}
\begin{proof}
We denote $G_i=\MT(A_i)$ the Mumford-Tate group of $A_i$, for $i=1,2$.
We first recall that if $G$ is a linear algebraic group over a field of characteristic $0$, and $V$ is a faithful finite dimensional representation of $G$, then every finite dimensional representation of $G$ can be constructed from the representation $V$ using the standard tannakian processes of forming tensor product, direct sums, subrepresentations, quotients and duals (see for example \cite{Wat}, 3.5). To sum up, we can say that the tannakian category $\Rep(G)$ of finite dimensional representations of $G$ is generated by $V$ (i.e. $\langle V\rangle \cong \Rep(G)$).
Let $V_i=\Hun(A_i,\Q)$ be the natural representation of $G_i$ for $i=1,2$.
Those two representations are faithful, so each one generates the corresponding tannakian category $\Rep(G_i)$.
Denote by $G$ the adjoint group $G_1^{ad}=G_2^{ad}$.
Let $V$ be a faithful representation of $G$. When we restrict this representation to $G_i$, we obtain a representation that is constructed from $V_i$ by the standard tannakian processes. But the category $\Rep(G_i)$ is equivalent to the category of absolute Hodge motives generated by the absolute Hodge motive $\hunah(A_i)$. So $V$ is the Hodge realization of two absolute Hodge motives: one constructed from $\hunah(A_1)$, denoted $M_1$ and the other constructed from $\hunah(A_2)$, denoted $M_2$. Replacing $F$ by a finite extension, we can suppose that all Hodge cycles on $A_i$ are defined over $F$. The Hodge realisations of $M_1$ and $M_2$ are equal by hypothesis. So $M_1$ and $M_2$ are isomorphic (possibly after a finite extension of $F$) as absolute Hodge motives because the Hodge functor from abelian absolute Hodge motives to Hodge structures is fully faithful over an algebraically closed field (see \cite{De1}, II, 6.25). This notably gives an isomorphism of the Galois representations associated to $M_1$ and $M_2$. We now suppose that $F$ is sufficiently big so that we have this isomorphism. We denote by $M$ the absolute Hodge motive $M_1$ ($\cong M_2$) over $F$.

Let $\p$ be a finite prime of $F$. We will now show that $M$ is potentially unramified at $\p$ if and only if $A_1$ has potentially good reduction at this prime.
Choose a prime $\ell$ of $\Q$ not equal to the residue characteristic of $\p$ and denote by $V_{1,\ell}$ the corresponding Galois representation (this is the $\ell$-adic realisation of the absolute Hodge motive $\hunah(A_1)$). By the monodromy theorem of Grothendieck \cite{SGA7}, Exp I, 3.6, the representation of the inertia group $\II(\bar{\p})$ of a prime $\bar{\p}$ of $\bar{F}$ over $\p$ is potentially unipotent on $V_{1,\ell}$. After a finite extension of $F$, we can suppose that this action is unipotent. It results from theorems of Deligne about absolute Hodge cycles \cite{De1}, Exp I, 2.9,2.11 that  after a finite extension of $F$, the galois representation on $V_{1,\ell}$ factors through the group of $\ql$-values of $G_1=\MT(A_1)$, i.e.
$$\II(\bar{v})\to G_1(\ql)\to \GL(V_{1,\ell}).$$
We now replace $F$ by a finite extension on which the inertia group acts unipotently and such that the Galois representation factors through the Mumford-Tate group.
Denote by $V_\ell$ the $\ell$-adic realisation of $M$.
Recall that the Neron-Ogg-Shafarevich criterion \cite{SeTa}, Th 1 says that $A_1$ has good reduction at $\p$ if and only if $V_{1,\ell}$ is unramified at $\p$. It remains to prove that $V_{1,\ell}$ is unramified at $\p$ if and only if $V_\ell$ is unramified at $\p$.

Suppose that $V_{1,\ell}$ is unramified at $\p$.
We know that $M$ is in the tannakian category $\langle\hunah(A_1)\rangle$, so the representation $V_\ell$ of the inertia group $\II(\bar{\p})$ is in the sub-tannakian category $\langle V_{1,\ell}\rangle$ of the category of $\ell$-adic representations of $\II(\bar{\p})$ generated by $V_{1,\ell}$. This implies that $V_\ell$ is unramified.
 
Suppose now that $V_\ell$ is unramified at $\p$. The natural map $\II(\bar{\p})\to G_1^{ad}(\ql)$ is trivial, so the inertia group acts on $V_{1,\ell}$ through the center of $G_1(\ql)$. Since $G_1$ is reductive, this center is of multiplicative type, so it does not contain unipotent points over $\ql$ (using our description \ref{unip-reductive} of unipotent points, this comes from \cite{SGA3II}, expXVII, 2.4). Thus the action of $\II(\bar{\p})$ on $V_{1,\ell}$ is trivial, i.e. $V_{1,\ell}$ is unramified.
\end{proof}

\begin{eprop}[Transposition statement]
\label{enonce-de-transposition}
Let $(G,X)$ be an adjoint Shimura datum with $G$ $\Q$-simple. If there exists a Deligne polymer $\mathscr{S}_0$ for $(G,X)$ such that every abelian variety associated to the triplet $(G,X,\mathscr{S}_0)$ has potentially good reduction at every discrete place of its field of definition, then for all Deligne polymer $\mathscr{S}$ for $(G,X)$, every abelian variety associated to the triplet $(G,X,\mathscr{S})$ has potentially good reduction at every discrete place of its field of definition.
\end{eprop}
\begin{proof}
Using \ref{constab-polymer}, we can associate to the Shimura datum $(G,X)$ and the Deligne polymer $\mathscr{S}_0$ a family
$\mathcal{A}_0 \to \Sh_{K_0}(G_0,X_0)_{E(G_0,X_0)}$ of abelian varieties over the canonical model of a Hodge type Shimura variety (here, $(G_0,X_0)$ is a Hodge type Shimura datum with adjoint $(G,X)$ and $K_0\subset G_0(\A_f)$ is a compact open subgroup). We do the same with another Deligne polymer $\mathscr{S}_1$ and obtain a family
$\mathcal{A}_1 \to \Sh_{K_1}(G_1,X_1)_{E(G_1,X_1)}$ of abelian varieties over the canonical model of a Hodge type Shimura variety ($K_1\subset G_1(\A_f)$ a compact open subgroup).

After possibly replacing the field $F$ by a finite extension that contains the composite field $E(G_0,X_0).E(G_1,X_1).E(G,X)$, and changing the compact open subgroups $K_0\subset G_0(\A_f)$, $K_1\subset G_1(\A_f)$ and $K\subset G(\A_f)$, we have a diagram
$$
\xymatrix{
\mathcal{A}_{0,F}\ar[d]		&		& \mathcal{A}_{1,F}\ar[d]\\
\Sh_{K_0}(G_0,X_0)_F\ar[dr]	&		& \Sh_{K_1}(G_1,X_1)_F\ar[dl]\\
			       	& \Sh_K(G,X)_F&\\
}
$$
Let $A_1:=\mathcal{A}_{1,x}$ be a fibre of the family $\mathcal{A}_{1,F}$ over a point $x_1\in \Sh_{K_1}(G_1,X_1)(F)$. The image of this point in $\Sh_{K}(G,X)(F)$ is denoted $x$. Take a point $x_0\in\Sh_{K_0}(G_0,X_0)(F)$ in the pre-image of the $G(\A_f)$-orbit of $x$ for the natural projection. Such a point exists because the projection
$$\Sh(G_0,X_0)(\C)\to \Sh(G,X)(\C)$$ is finite (see \cite{De4}, 2.1.2, 2.1.8) and surjective on a connected component, and the action of $G(\A_f)$ on the set of connected components of $\Sh(G,X)(\C)$ is transitive (see \cite{De4}, 2.1.3).
The fibre of the family $\mathcal{A}_0$ over $x_0$ is an abelian variety $A_0$.
All the Hodge structures of the points of the orbit $G(\A_f).x_0$ are equal because the $G(\A_f)$-action doesn't change the Hodge structure.
This shows that the abelian varieties $A_0$ and $A_1$ verify the hypothesis of 
\ref{adjointtip} so $A_0$ has potentially good reduction at a place of $F$ if and only if $A_1$ has potentially good reduction at this place. By construction, $A_0$ is an abelian variety associated to the triplet $(G,X,\mathscr{S}_0)$. The hypothesis implies that $A_0$ has potentially good reduction at every discrete place of $F$. This implies that $A_1$ has potentially good reduction at every discrete place of $F$.
\end{proof}

\subsection{Explicit examples}
\label{explicit-examples}
\begin{edefin}
\label{definpotperftenstord}
Let $(G,X)$ be a Shimura datum of classical adjoint type with $G$ simple. We will say that \emph{$(G,X)$ is of potential perfectly tens-twisted type} if there exists a Deligne polymer $\mathscr{S}_1$ for $(G,X)$ such that the representation $G_1\to \GL(V)$ given by the corresponding Siegel embedding
$$(G_1,X_1)\to (\GSp(V),S^\pm)$$
is perfectly tens-twisted over $\ql$ for some prime $\ell$ of $\Q$.
\end{edefin}
\begin{eprop}
\label{proppotperftens}
Let $(G,X)$ be a Shimura datum of potential perfectly tens-twisted type.
Any abelian variety associated to any Deligne polymer $\mathscr{S}$ for $(G,X)$ has potentially good reduction at every discrete place of its field of definition.
\end{eprop}
\begin{proof}
By the transposition statement \ref{enonce-de-transposition}, we know that it is enough to show the result for one Deligne polymer. Therefore, we will show this for the Deligne polymer $\mathscr{S}_1$ given by hypothesis. The result on potentially good reduction comes directly from our combinatorial criterion \ref{combincriteria} and from the constructive result \ref{constructive-result}.
\end{proof}

\begin{edefin}
Let $(G,X)$ be a Shimura datum of classical adjoint type (i.e. $G$ is $\Q$-simple of type $A$,$B$,$C$ or $D$) and suppose that $G_\R$ has at least two simple factors and that only one factor of $G_\R$ is non-compact. Such a Shimura datum will be called a \emph{potential Mumford type Shimura datum}.
\end{edefin}

Let $(G,X)$ be a potential Mumford type Shimura datum. There exists a Deligne polymer $\mathscr{S}_1$ for $(G,X)$ containing at least one part $T\subset \mathscr{P}(D)$ that projects surjectively on the set $I$ of simple factors of $G_\R$.

The construction will be done in several simple steps. 
First, denote by $D_v$ the part of the Dynkin diagram corresponding to the simple non compact factor $G_v$ of $G_\R$. Let $s_v$ be the vertex of $D_v$ corresponding to $X$ (in the notations of \cite{De4}, 1.2.6). We denote $s_{T,v}$ an underlined vertex of $D_v$ in the table \cite{De4}, 1.3.9 for $(D_v,s_v)$.
For each compact factor $G_w$ of $G_\R$, we have an isomorphism $D_w\cong D_v$ and we denote $s_{T,w}$ the vertex corresponding to $s_{T,v}$ by this isomorphism. The set $\{s_{T,w}\,|\,w\in I\}$ of all those vertices gives a part $T\subset\mathscr{P}(D)$ of the Dynkin diagram of $G$ that verifies the
condition of the definition of a Deligne polymer \ref{defindelignepolymer}. Indeed, $T\cap D_{nc}=\{s_{T,v}\}$ is reduced to one point and $s_{T,v}$ is an underlined vertex of $D_v$ in the table \cite{De4}, 1.3.9 for $(D_v,s_v)$.
Let $\mathscr{S}_1$ be the $\Gal(\qb/\Q)$-orbit of $T$ in $\mathscr{P}(D)$. Every $T'$ in $\mathscr{S}_1$ verifies the condition of the definition of a Deligne polymer \ref{defindelignepolymer}. Indeed, $T'\cap D_{nc}$ is empty or reduced to one point $s_T\in D_v$ ($v\in I_{nc}$) and $s_T$ is in the galois orbit of some $s_{T,w}$ so $s_T=s_{T,v}$ or $s_T=s_{T,v}^*$ where $*$ denotes the opposition involution of the Dynkin diagram $D_v$. Since the underlined vertex of \cite{De4}, 1.3.9 are stable by the opposition involution, $s_T$ is an underlined vertex for $(D_v,s_v)$ in \cite{De4}, 1.3.9. This shows that $\mathscr{S}_1$ verifies the condition of the definition of a Deligne polymer \ref{defindelignepolymer}.
Moreover $\mathscr{S}_1$ contains a non-empty part and is Galois stable by construction. It is covering ($\cup_{T'\in\mathscr{S}_1}\pi(T')=I$ for $\pi\colon D\to I$ the natural projection) because $G$ is $\Q$-simple so $\Gal(\qb/\Q)$ acts transitively on $I$.
So $\mathscr{S}_1$ is a Deligne polymer and $T$ projects surjectively on the set $I$ of simple factors of $G_\R$.
\begin{eprop}
\label{proppotmum}
Let $(G,X)$ be a potential Mumford type Shimura datum.
Every abelian variety associated to every Deligne polymer $\mathscr{S}$ for $(G,X)$ has potentially good reduction at every discrete place of its field of definition.
\end{eprop}
\begin{proof}
We will show that the Deligne polymer constructed in the precedent paragraph verifies the condition given in the definition \ref{definpotperftenstord}, so gives $(G,X)$ the structure of a potential perfectly tens-twisted type Shimura datum. It suffices to show that there exists a prime $\ell$ of $\Q$ such that the representation $G_1\to\GL(V)$ given by a Siegel embedding $(G_1,X_1)\to (\GSp(V),S^\pm)$ associated to $\mathscr{S}_1$ is perfectly tens-twisted over $\ql$. This was already done in the study of Mumford type reprentations in \ref{cormumford}. After that, we apply the result given in \ref{proppotperftens} to the potential perfectly tens-twisted Shimura datum $(G,X)$ to conclude the proof.
\end{proof}

\begin{ecor}
Let $A$ be a simple abelian variety over a number field with non-commutative Mumford-Tate group and suppose that $\MT(A)^{ad}$ is $\Q$-simple,
$\MT(A)^{ad}_\R$ has at least two simple factors and only one non-compact factor. Then $A$ has potentially good reduction at every discrete place of its field of definition.
\end{ecor}
\begin{proof}
$(\MT(A),X_A)^{ad}$ is a potential Mumford type Shimura datum so the conjugation of proposition \ref{proppotmum} and of \ref{exmabfibre} is sufficient to obtain the result.
\end{proof}

\begin{edefin}
Let $F$ be a number field.
Let $\norm\colon F\to \Q$ be the usual norm defined by $\norm(f)=\det(x\mapsto f.x)$.
Following Ferrand\cite{Ferrand}, we define the norm or corestriction of $D$ from $F$ to $\Q$ as a couple
$(\mathrm{N}_{F/\Q}D,\nu_D)$ of
a central simple algebra $\mathrm{N}_{F/\Q}D$ and a multiplicative functorial map $\nu_D\colon D\to \mathrm{N}_{F/\Q}D$
between the Algebra-schemes corresponding to those algebras that verifies
$$\nu_A(fx)=\norm_A(f)\nu_A(x)$$
for all $f\in A$ and $x\in D_A$ with $A$ a commutative $\Q$-algebra and that is universal for this property.
\end{edefin}

\begin{eexm}
Let $F$ be a totally real number field of degree equal to $2n+1$, $n>0$ over $\Q$.
Let $D$ be a quaternion algebra over $F$. Suppose that $D$ decomposes over $\R$ as a product
$$D_\R\cong \mathrm{M}_{2,\R}\times \mathbb{H}^{2n}$$
of a matrix algebra and a product of quaternion algebras.
Let $G_2=\Res_{F/\Q}D^\times$, $G=G_2^{ad}$ and $h_2\colon \mathbb{S}\to G_{2,\R}$ be the trivial morphism on quaternionic parts and the classical morphism
$$
\begin{array}{cccc}
h_0: & \mathbb{S}	& \to 	& \GL_{2,\R}\\
     & a+ib             &\mapsto& 	\left(\begin{array}{cc}
   					a  & b\\
					-b & a
	   				\end{array}\right)
\end{array}
$$
on the matrix part. Let $h\colon \mathbb{S}\to G_\R$ be the projection of $h_2$. We suppose that $\mathrm{N}_{F/\Q}D$ splits, i.e. is isomorphic to
$\mathrm{M}_{2^{2n+1},\Q}$. By definition, the morphism $\nu_D$ induces a canonical morphism of group schemes
$$\nu_D^\times\colon \Res_{F/\Q}D^\times\to\GL(V),$$
with $V=\Q^{2^{2n+1}}$.
Denote $G_1$ the image of this morphism and $h_1\colon \mathbb{S}\to G_{1,\R}$ the projection of $h_2$. Denote $X$ (resp. $X_1$, resp. $X_2$) the conjugacy class of morphisms corresponding to $h$ (resp. $h_1$, resp. $h_2$). Then $(G,X)$ and $(G_1,X_1)$ are Shimura data and since the kernel of $G\to\GL(V)$ is central, $(G_1,X_1)^{ad}=(G,X)$. The natural involution $\sigma_D$ of $D$ gives by functoriality an involution on $\mathrm{N}_{F/\Q}D$. This involution is symplectic because it is equal over $\qb$ to a tensor product of $2n+1$ involutions of symplectic type (the classical symplectic involution on $\mathrm{M}_{2,\qb}$). So we have a diagram
$$(G,X)\gets(G_1,X_1)\hookrightarrow (\GSp(V),S^\pm).$$
The projection $\pi\colon D\to I$ of the Dynkin diagram $D$ of $G$ to the set $I$ of simple factors of $G$ is a bijection because each factor of $G_\qb$ is of $A_1$ type.
The polymer of the above diagram contains only one part:
$$\mathscr{S}_1=\mathscr{S}(V)=\{D\}.$$
We thus have constructed a Shimura datum $(G,X)$ that is of potential Mumford type. We know from \ref{proppotmum} that the abelian varieties that are the fibres over number fields of the family associated to the prescribed Siegel embedding of $(G_1,X_1)$ have potentially good reduction at every finite prime of their field of definition.

If we choose another polymer, for example the $\Gal(\qb/\Q)$-orbit $\mathscr{S}_{\langle E\rangle}$ of any non empty subset $E$ of $D$, we also know from \ref{proppotmum} that any abelian variety associated to $(G,X,\mathscr{S}_{\langle E\rangle})$ has potentially good reduction at every discrete place of its field of definition. For example if $E=\{d_0\}$, we obtain a Shimura variety that is not far from a {\PEL} Shimura variety (see \cite{Reimann} for a detailled study of this case).
\end{eexm}

\begin{erem}
The result of the last example can also be proved using the transposition result \ref{enonce-de-transposition}, the valuative criterion of properness and a properness result of Reimann \cite{Reimann}, 2.14 that uses methods similar to Morita's method for Shimura varieties that are allmost {\PEL} Shimura varieties (weak polarization).
It is possible to construct, using Deligne polymers, a similar example with a group $G$ of $B_n$ type. In this case, we know from \ref{Bn-paspel} that there is no hope to reduce the question to some {\PEL} Shimura variety. We preferred to give the simplest example of a group of type $A_1$, but our method really gives new results that can not be reduced to {\PEL} results using \ref{enonce-de-transposition}.
\end{erem}

\begin{edefin}
\label{defin-pot-well-twisted}
Let $(G,X)$ be a Shimura datum of classical adjoint type with $G$ simple. We will say that \emph{$(G,X)$ is of potential well twisted cyclic type} if the two following conditions are verified:
\begin{itemize}
\item The action of $\Gal(\qb/\Q)$ on the set $I$ of simple factors of $G_\qb$ is cyclic.
\item There exists a Deligne polymer $\mathscr{S}_1$ for $(G,X)$ containing at least one part $T\subset \mathscr{P}(D)$ of the Dynkin diagram $D$ of $G$ with $\Card(T)>1$.
\end{itemize}
\end{edefin}
\begin{eprop}
\label{proppotcyclic}
Let $(G,X)$ be a Shimura datum of potential well twisted cyclic type.
Every abelian variety associated to every Deligne polymer $\mathscr{S}$ for $(G,X)$ has potentially good reduction at every discrete place of its field of definition.
\end{eprop}
\begin{proof}
We will show that the Deligne polymer $\mathscr{S}_1$ given by definition  verifies the condition given in \ref{definpotperftenstord}, so gives $(G,X)$ the structure of a potential perfectly tens-twisted type Shimura datum.
It suffices to show that there exists a prime $\ell$ of $\Q$ such that the representation $G_1\to\GL(V)$ given by a Siegel embedding $(G_1,X_1)\to (\GSp(V),S^\pm)$ associated to $\mathscr{S}_1$ is perfectly tens-twisted over $\ql$. This was already done in \ref{corcyclic}.
After that, we apply the result given in \ref{proppotperftens} to the potential perfectly tens-twisted Shimura datum $(G,X)$ to conclude the proof.
\end{proof}

\begin{eexm}
Let $F$ be a totally real number field with a fixed isomorphism $\Gal(F/\Q)\to\Z/8\Z$.
Let $D$ be a quaternion algebra over $F$. We know that $D$ decomposes
over $\R$ as a product
$$D_\R\cong \mathrm{M}_{2,\R}^m\times \prod \mathbb{H}^{n-m}$$
of $m$ matrix algebras and $m-n$ quaternion algebras, for $0\leq m\leq 8$.
Let $G_2=\Res_{F/\Q}D^\times$, $G=G_2^{ad}$ and $h_2\colon \mathbb{S}\to G_{2,\R}$ be
the trivial morphism on quaternionic parts and the classical morphism
$$
\begin{array}{cccc}
h_0: & \mathbb{S}	& \to 	& \GL_{2,\R}\\
     & a+ib             &\mapsto& 	\left(\begin{array}{cc}
   					a  & b\\
					-b & a
	   				\end{array}\right)
\end{array}
$$
on the matrix part. Let $h\colon \mathbb{S}\to G_\R$ be the projection of $h_2$. Let $X$ and $X_2$ be the corresponding conjugacy classes of morphisms.
The Shimura datum $(G,X)$ is of adjoint type. The projection $\pi\colon D\to I$ from the Dynkin diagram to the set of simple factors of $G_\qb$ is a bijection because each factor of $G_\qb$ is of $A_1$ type. We choose a bijection $I\to \{1,\dots,8\}$ compatible with the isomorphism $\Gal(F/\Q)\to\Z/8\Z$ and we will now use this identification without mentioning it. Denote by $I_{nc}$ the subset of $I$ corresponding to noncompact factors of $G_\R$. Using that $\pi\colon D\to I$ is a bijection, we deduce that the exhaustive classification given by Addington \cite{Add}, 6.3 of all of ``her polymers'' in this case gives the classification of all Deligne polymers for $(G,X)$.
We can see for example for $I_{nc}=\{1,5\}$ that the Deligne polymer
$$\mathcal{S}_1=\{1234,2345,3456,4567,5678,6781,7812,8123\}$$
verifies the second hypothesis of \ref{defin-pot-well-twisted}. So $(G,X)$ is of potential well twisted cyclic type and any abelian variety associated to any Deligne polymer for $(G,X)$ (their list is given in \cite{Add}, 6.3) has potentially good reduction at every finite prime of its field of definition.
\end{eexm}

\begin{edefin}
Let $(G,X)$ be a Shimura datum of classical adjoint type with $G$ simple. We will say that \emph{$(G,X)$ is of potential {\PEL} type} if there exists a Deligne polymer $\mathscr{S}_1$ for $(G,X)$ such that a corresponding Shimura datum $(G_1,X_1)$ is of {\PEL} type.
\end{edefin}
\begin{eprop}
\label{proppotPEL}
Let $(G,X)$ be a Shimura datum of potential {\PEL} type such that there exists a place $v$ of $F$ such that $G_v$ is anisotropic.
Every abelian variety associated to every Deligne polymer $\mathscr{S}$ for $(G,X)$ has potentially good reduction at every discrete place of its field of definition.
\end{eprop}
\begin{proof}
The results of Morita \cite{Morita} about his conjecture in the {\PEL} case, with the additional $\lambda$-adic result we gave in \cite{TheseFred}, Part 2 gives that every abelian variety associated to the Deligne polymer $\mathscr{S}_1$ (given by definition) for $(G,X)$ has potentially good reduction at every discrete place of its field of definition. The transposition statement \ref{enonce-de-transposition} permits to conclude.
\end{proof}

Institut Math\'ematique de Rennes, universit\'e de Rennes 1,
Campus de Beaulieu, 35042 Rennes.\\
e-mail: frederic.paugam@math.univ-rennes1.fr

\bibliographystyle{alpha}
\bibliography{$HOME/travail/fred.bib}

\end{document}